\documentclass[journal]{IEEEtran}

\ifCLASSINFOpdf
   \usepackage[pdftex]{graphicx}

\else

   \usepackage[dvips]{graphicx}

\usepackage{amsmath,amssymb}
\usepackage{graphicx}
\usepackage{cite}
\usepackage{float}
\usepackage{bm}
\usepackage{caption}
\renewcommand{\theequation}{\arabic{section}.\arabic{equation}}

\newtheorem{theoremx}{Theorem}[section]
\newtheorem{lemmax}{Lemma}[section]

\newtheorem{remark}{Remark}[section]
\newtheorem{assumption}{Assumption}[section]

\usepackage{chngcntr}
\counterwithout{equation}{section}

\ifCLASSINFOpdf
\else
\fi
\begin{document}

\title{Distributed stochastic subgradient-free algorithm for Nash equilibrium seeking in two-network zero-sum games}
\author{Dandan Yue, Ziyang Meng
\thanks{This work has been supported in part by National Natural Science
Foundation of China under Grants 61803222, 61873140 and 61833009, in part by the Independent Research Program of Tsinghua University under Grant 2018Z05JDX002.
}
\thanks{D. Yue and Z. Meng are with Department of Precision Instrument,
Tsinghua University, Beijing 100084, China. Corresponding author: Ziyang Meng. ddyue@mail.tsinghua.edu.cn (D. Yue), ziyangmeng@mail.tsinghua.edu.cn (Z. Meng). }}

\maketitle

\begin{abstract}
This paper investigates the distributed Nash equilibrium seeking problem for two-network zero-sum games with set constraints, where the two networks have the opposite nonsmooth cost functions. The interaction of the agents in each network is characterized by an unbalanced directed graph. We are particularly interested in the case that the local cost function of each agent in the two networks is unknown in this paper. We first construct the stochastic subgradient-free two-variable oracles based on the measurements of the local cost functions. Instead of using subgradients  of the local cost functions, the subgradient-free two-variable oracles are employed to design the distributed algorithm for the agents to search the Nash equilibrium. Under the strong connectivity assumption, it is shown that the proposed algorithm guarantees that the states of all the agents converge almost surely to the neighborhood of the Nash equilibrium, where the scale of the neighborhood can be arbitrarily small by selecting suitable smooth parameter in the oracles. Numerical simulations are finally given to  validate the theoretical results.
\end{abstract}

%????????????two networks engaged in the zero-sum games consist of the strategic agents
%
% and weight-unbalanced
%
%, where the neighbors' out-degree of each agent is not required


\IEEEpeerreviewmaketitle

%he weight matrices are only row-stochastic
%
%The algorithm does not require each agent to know their neighbors' out-degree.
%
%problem, novel, setting, other paper, results
%
% The cost function of each network is the sum of the local cost functions associated with the agents in the respective network. Each of the two networks is equipped with a directed and weight-unbalanced graph.

\section{Introduction}

Zero-sum game problem has arisen in engineering application domains such as adversarial estimation of sensor networks and power allocation of channels in presence of jammers \cite{Vamvou12,aram12}. Due to the theoretical significances and broad applications, the study on how to search the Nash equilibrim has been one of the main focuses of the works on zero-sum games \cite{Mazum19}. Various approaches have been developed to find the Nash equilibrium in zero-sum games, including fictitious play \cite{Shamm04}, best response approach  \cite{hofb06}, and subgradient-based approach \cite{Nedic09}. The above works \cite{Shamm04,hofb06,Nedic09} solve the problems in the centralized settings and consider that all the strategic agents are adversarial. Recently, the study on zero-sum games in the networked settings has received much attention. In particular, Gharesifard and Cortes \cite{Ghare13} focused on the two-network zero-sum games, where the two networks, viewed as two players, were adversarial while the strategic agents in the same network were collaborative.  Based on the local information exchange, a distributed method that synthesizes the consensus dynamics and the saddle point dynamics was designed to seek the Nash equilibrium. Following the distributed framework proposed in \cite{Ghare13}, many extended works have been presented to  solve the Nash equilibrium seeking problem in two-network zero-sum games. In \cite{yang18}, the  distributed projection gradient method was proposed for the distributed constrained minimax problem, which can be viewed as the constrained two-network zero-sum game problem from the perspective of game theory. Then, the authors of \cite{shic19} focused on the case of the sequential communication and an incremental algorithm was proposed to enable the agents in different networks to update their strategies asynchronously. Note that the graphs of the two networks are  fixed and undirected in \cite{yang18} and \cite{shic19}, and the authors of  \cite{Lou16} developed a distributed subgradient algorithm that works for the switching and unbalanced directed graphs.

The implementation of the aforementioned distributed algorithms requires the gradient or subgradient information of the agents' cost functions. However, such information may be not available when the objective/cost functions are unknown, or it is costly to obtain in practical applications. In this paper, we focus on the two-network zero-sum
game problem where the local cost function of each agent is unknown. In fact, the gradient-free schemes have been used to solve the problems of optimizations and games in presence of the unknown cost functions. In particular, the authors of  \cite{Nester17} and  \cite{pang19}  introduced the stochastic gradient-free methods to solve the optimization problems. Ref. \cite{Nester17}  focused on the centralized method and ref. \cite{pang19} extended it to a distributed method over a directed graph.  A stochastic gradient-free algorithm generated by the stochastic differential inclusion was proposed for the noncooperative games in \cite{poveda15}. In \cite{zhu16}, the finite differences were used to approximate the gradients. The evolution strategy, which is developed based on the stochastic gradient approximator, was used in \cite{Alduj19} to solve the (centralized) minimax (or called two-player zero-sum game) problem. Note that the cost functions involve only one vector in \cite{Nester17} and \cite{pang19}  while each agent's local cost function depends on two vectors corresponding to two networks in two-network zero-sum games. Therefore, the stochastic gradient-free methods  presented in \cite{Nester17} and \cite{pang19} can not be directly used for the two-network zero-sum games. Although the methods developed in \cite{poveda15} and \cite{zhu16} solve games (including zero-sum games), they require the objective functions to be quadratic or smooth. The problem studied in \cite{Alduj19} is similar to that considered in this paper, however, it is solved only in a centralized setting. Moreover, even in the networked settings, the methods presented in \cite{poveda15,zhu16,Alduj19} are not applicable to the two-network zero-sum games since the agents (in the same network) can not reach agreement about their estimations of the corresponding vector involved in the network cost function. Motivated by the above discussions, this paper aims at developing the distributed subgradient-free algorithm for the two-network zero-sum games.

The main contributions of this paper are three-fold. First, a distributed stochastic subgradient-free algorithm is proposed to achieve the Nash equilibrium seeking of two-network zero-sum games. Unlike existing works on  the  two-network zero-sum games \cite{Ghare13,yang18,shic19,Lou16}, no exact mathematical formulations but only the measurements of the local cost functions are needed herein. Second, different from the strategies adopted in \cite{poveda15,zhu16,Alduj19}, the proposed algorithm solves  the problem with the general nonsmooth cost functions in a distributed manner. Third, from a technical perspective, the two-variable oracles are properly designed in this paper for the implementation of the proposed algorithm. This can be viewed as an extension of one-variable oracles given in refs. \cite{Nester17,pang19}.

The structure of this paper is as follows. The preliminaries are provided and the problem is introduced in Section II. The proposed algorithm is presented in Section III. The almost sure convergence of the proposed algorithm is investigated in Section IV. The numerical simulations are presented in Section V and the conclusion is given in Section VI.

Notations:  For a differentiable function $f$, let $\nabla_{x} f(\cdot,y)$ and $\nabla_{y} f(x,\cdot)$ denote the gradients of $f(\cdot,y)$ at $x$ and $f(x,\cdot)$ at $y$, respectively.  For a locally Lipschitz and convex but not  necessarily differentiable function $f$, let ${{\partial _x}{f}(\cdot,y)}$ and ${{\partial _y}{f}(x,\cdot)}$ denote the subgradients (or generalized gradients) of $f(\cdot,y)$ at $x$ and $f(x,\cdot)$ at $y$, respectively. The following are from refs. \cite{Lukacs14,Rotar06,polyak87}. For the random variables $a$ distributed on the set ${\cal U} \subseteq {\mathbb R}^{n_1}$ and $b$ distributed on the set ${\cal W} \subseteq {\mathbb R}^{n_2}$, $\Xi [a|b = \tilde b] = \int_{\cal U} \tilde a\psi  (\left. {\tilde a} \right|\tilde b)d\tilde a$ denotes the conditional expectation of $a$ given $b=\tilde b$ and $\Xi[a] \buildrel \Delta \over = \Xi\big[\Xi [a|b] \big]$ denotes the total expectation of $a$, where $\psi  (\left. {\tilde a} \right|\tilde b)$ is the conditional probability density of $a$ given $b=\tilde b$, and $\tilde a$ and $\tilde b$ are real vectors in ${\cal U}$ and ${\cal W}$, respectively.  A random event happens ``almost surely" means that it happens ``with probability one". We will use these two terminologies alternatively, and abbreviate ``almost surely" by a.s. in this paper. ${\rm I}_i$ is a vector where the $i$-th element is $1$ and other elements are $0$. $0_n$ denotes an $n$-dimension vector consisting of $0$. $\left|  \cdot  \right|$ and $\left\|  \cdot  \right\|$ refer to the absolute value and the Euclidean norm of a scalar and a vector, respectively.

% The collaborative agents were in charge of computing the vector, which corresponded to their own network, involved in the network cost function.

%we must give the papers that use the gradient-free method in the game and use the balanced graph.
%
%is nontrival from using the gradient technique from the one-variable optimization to the game???????

 %In practice, an entire objective is usually achieved via the cooperation of the individuals. For example, the multiple channels make decisions on the signal powers to maximize the secrecy capacity of the whole communication system in presence of the noise interference \cite{boyd04}.
%
%in fact, there are some works that for distributed NE of the two-player zero-sum game and gradient-free, they use such as the neural networks and the adaptive dynamic programming to...see the HJI equation. but ???????see the lack of this in zero-order TAC paper?????????
%
%When we describe the existing works, we should pay attention that stochastic-matrix is balance condition!!!!!!
%
%When we describe the contibution about the only using row-stochastic matrix, we should exphazie on the graph since this sounds more applicable and valuable, and say that as a result ... the balance matrix is the theory espect induce by the graph(see Introduction of TAC-Linear convergence )

\section{Preliminaries and Problem Statement}

\subsection{Graph theory}

Consider a directed graph ${\cal G}=({\cal V},{\cal E})$, with a node set ${\cal V}$ and an edge set ${\cal E} \subseteq {\cal V} \times {\cal V}$. For any $s,i \in   {\cal V}$, $(s,i) \in {\cal E}$ means that node $s$ sends its information to node $i$, and $s$ is called the neighbor (i.e., in-neighbor) of $i$. The adjacency matrix of ${\cal G}$ is defined by $A=[{a_{is}}]$, where ${{a_{is}}}>0$ for $s \in {\cal N}_{i}$ and ${{a_{is}}}=0$ for $s \in {\cal V} / {\cal N}_{i}$, with $ {\cal N}_{i}$ denoting the set of neighbors (i.e., in-neighbors) of $i$ in ${\cal V}$ (including itself). $A$ is a row-stochastic matrix if $\sum\nolimits_{s \in {{\cal V}}} {a_{is}} =1$ for any $i \in {{\cal V}}$.  A path from $s$ to $i$ is a sequence of edges $(s,s_1), \ldots,(s_n,i)$. A directed graph is strongly connected if there exists a path from any node $s$ to any node $i$ for $s \ne i$.

%The adjacency matrix of ${\cal G}$ is $A=[a_{ij}]$, where $a_{ij}>0$ if $(j,i) \in {\cal E}$, and $a_{ij}=0$ otherwise.

%The set of the neighbors (in-neighbors) of $i$ is denoted as $ {\cal N}_i$.

\subsection{Convex analysis}

${P_{\cal X} }:  \mathbb{R}^n \to {\cal X} $ is the projection operator onto a closed convex set ${\cal X} \subseteq \mathbb{R}^n$ if for any vector $ x \in \mathbb{R}^n $, $P_{\cal X} (x)$ satisfies $P_{\cal X} (x) \in {\cal X} $ and $\|x- P_{\cal X} (x)\|=\inf_{z \in {\cal X} }\| x-z\|$. A locally Lipschitz function $f(\cdot):\mathbb{R}^n \rightarrow \mathbb{R}$ has the generalized gradient $\partial_x f$. If $f$ is convex, $\partial_x f$ coincides with the subdifferential of $f$ at $x$, where the subdifferential is the set of all subgradients, i.e., $\partial_x f=\{ z \in \mathbb{R}^n: f(y)-f(x) \ge z^T(y-x), \forall y \in \mathbb{R}^n\}$. A function $f$ is concave if $-f$ is convex. $f(x,y):\mathbb{R}^{n_1} \times \mathbb{R}^{n_2} \rightarrow \mathbb{R}$ is (strictly) convex-concave if $f(x, y)$ is (strictly) convex with respect to $x$ for any fixed $y \in \mathbb{R}^{n_2} $ and is (strictly) concave with respect to $y$ for any fixed $x \in \mathbb{R}^{n_1} $.

\subsection{Saddle point, zero-sum game, and Nash equilibrium}

A pair $(x^*,y^*)$ is a saddle point of $F(x,y)$ on ${\cal X} \times {\cal Y} \subseteq {\mathbb R}^{n_1} \times {\mathbb R}^{n_2}$ if $F(x^*,y) \le F(x^*,y^*)  \le F(x,y^*)$ for any $(x,y) \in {\cal X} \times {\cal Y}$. Suppose that ${\cal X} \subseteq {\mathbb R}^{n_1}$ and ${\cal Y}  \subseteq {\mathbb R}^{n_2}$ are two closed convex sets and $F(x,y)$ is a strict convex-concave function. Then, $F$ has a unique saddle point  \cite{Ghare13,Lou16}.

Suppose that the set of players in a game is denoted by ${\cal V}= \{ 1,2, \ldots ,m\} $ and the cost function of player $i \in {\cal V}$ is $F_i \buildrel \Delta \over = {F_i}({z_i},{z_{ - i}})$, where $z_i$ is the strategy of player $i$, ${z_{ - i}}=(z_1,  \ldots, z_{i-1},z_{i+1}  \ldots, z_m)$, and $m$ is the number of players. The game is called a zero-sum game if $\sum\nolimits_{i \in {\cal V}} {F_i}({z_i},{z_{ - i}})=0$. Let $ {\cal Z}_i$ be the convex strategy set of player $i$. A profile $z^*=(z^*_i,z^*_{-i})$ is a Nash equilibrium if $F_i({z^*_i},{z^*_{ - i}}) \le F_i({z_i},{z^*_{ - i}})$ holds for all $z_i \in {\cal Z}_i$ and $i \in {\cal V}$. For a two-player zero-sum game, a point pair $(z_1^*,z_2^*)$ is the Nash equilibrium if and only if it is the saddle point of the cost function $F_1$ \cite{Lou16}.

%For a two-player zero-sum game, let the Nash equilibrium $(z_1^*,z_2^*)$ be the Nash equilibriuma. Clearly, $(z_1^*,z_2^*)  \in {\cal Z}_1 \times {\cal Z}_2$, and it holds that $u_1(z_1^*,z_2^*) \le u_1({z_1},z_2^*) $ for any $z_1 \in {\cal Z}_1$ and $-u_1(z_1^*,z_2^*)=u_2(z^*_1,z_2^*) \le u_2({z^*_1},z_2)= -u_1(z_1^*,z_2)$ for any $z_2 \in {\cal Z}_2$. These two relations are equivalent  $u_1(z_1^*,z_2) \le u_1(z_1^*,z_2^*) \le u_1({z_1},z_2^*) $ for any $(z_1,z_2) \in {\cal Z}_1 \times {\cal Z}_2$.

%\section{Problem statement and algorithm description}

\subsection{Problem statement}

Consider a network with $m_1+m_2$ agents and consisting of two subnetworks. The agents in the two subnetworks are labeled as $1, 2,\ldots ,{m_1}$ and $1, 2,\ldots ,{m_2}$, respectively. Each agent in the two subnetworks knows the numbers $m_1$ and $m_2$, respectively. Denote ${{\cal V}_1} = \{ 1,2, \ldots ,{m_1}\}$ and ${{\cal V}_2} = \{ 1,2, \ldots ,{m_2}\} $. Let $f_{1i}(x,y):{{\mathbb R}^{{n_1}}} \times {{\mathbb R}^{{n_2}}} \to {\mathbb R}$ be the (possibly nonsmooth) local cost function of each agent $i \in {\cal V}_1$. The global cost function of the first subnetwork, i.e.,  $F(x,y)$,  is the sum of all these local cost functions:
\begin{equation}
F(x,y)=\sum\nolimits_{i \in {{\cal V}_1}} {{f_{1i}}(x,y)}. \label{1}
\end{equation}
Similarly, let $f_{2j}(x,y):{{\mathbb R}^{{n_1}}} \times {{\mathbb R}^{{n_2}}} \to {\mathbb R}$ be the (possibly nonsmooth) local cost function of each agent $j \in {\cal V}_2$. The global cost function of the second subnetwork, i.e., $\tilde F(x,y)$,  is the sum of all these local cost functions:
\begin{equation}
\tilde F(x,y)= - \sum\nolimits_{j \in {{\cal V}_2}} {{f_{2j}}(x,y)}. \label{1-1-1}
\end{equation}
We consider that the cost functions of these two subnetworks satisfy
\begin{equation}
F(x,y) + \tilde F(x,y)=0 \label{1-1-1-1}
\end{equation}
for any $(x,y) \in {\cal X} \times {\cal Y}$. This implies that $F(x,y)=\sum\nolimits_{j \in {{\cal V}_2}} {{f_{2j}}(x,y)}$. Note that the two subnetworks act as two players and are engaged in a two-player zero-sum game. We call such a game as a two-network zero-sum game  \cite{Ghare13,Lou16}. The strategy (constraint) set of this game is ${\cal X} \times {\cal Y}$, where  ${\cal X} \subseteq {{\mathbb R}^{{n_1}}}$ and ${\cal Y} \subseteq {{\mathbb R}^{{n_2}}}$ are the closed and convex sets. In this paper, the exact mathematical expressions of all the local cost functions $f_{1i}$ and $f_{2j}$ are not available while each agent $i \in {\cal V}_1$ and $j \in {\cal V}_2$ have access to the measurements of $f_{1i}$ and $f_{2j}$, respectively. To be more specific,  given any pair $(x,y) \in {{\mathbb R}^{{n_1}}} \times {{\mathbb R}^{{n_2}}}$, the value of $f_{1i}(x,y)$ can be obtained by agent $i \in {\cal V}_1$ and the value of $f_{2j}(x,y)$ can be obtained by agent $j \in {\cal V}_2$. Also, denote ${\cal Z}^*={\cal X}^* \times {\cal Y}^*$ as the nonempty set composed of the saddle point  $(x^*,y^*)$ of $F(x,y)$ on ${\cal X} \times {\cal Y}$.

\begin{figure}[t]
\centering
\includegraphics[height=2.2cm ,width=4.5cm]{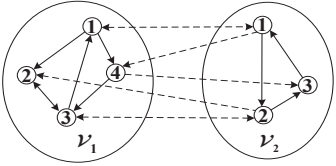}
\caption{An illustrative example for the graph in the whole network. The solid lines denote the communication links within ${\cal V}_1$ and ${\cal V}_2$, and the dotted lines denote the communication links between ${\cal V}_1$ and ${\cal V}_2$.}
\label{fig1}
\end{figure}

For the setup of graphs, the agents share the information with their neighbors  in the same subnetwork, meanwhile, the two subnetworks have access to the information about each other. We describe the communications among the agents in the whole network by a directed bipartite graph ${\cal G}$ consisting of adjoint node sets ${\cal V}_1$ and ${\cal V}_2$, where each node has at least one neighbor (not including itself). The communications within each subnetwork are described by directed graphs ${\cal G}_1$ and ${\cal G}_2$, respectively. Each node in ${\cal V}_1$ (or ${\cal V}_2$) has at least one neighbor from ${\cal V}_2$ (or ${\cal V}_1$).  An illustrative example for the graph of the whole network is given in Fig. {\ref{fig1}}. In what follows, we let  $A_1$ and $A_2$ be the adjacency matrices associated with  ${\cal G}_1$ and ${\cal G}_2$, which are row-stochastic.

The objective of this paper is to design a distributed stochastic subgradient-free algorithm such that the Nash equilibrium seeking is achieved for the considered two-network zero-sum games.

\section{Distributed stochastic subgradient-free algorithm}

\subsection{Two-variable Gaussian approximations of cost functions}

Before presenting the proposed distributed stochastic subgradient-free algorithm, we first define the smoothed versions of (\ref{1}) and (\ref{1-1-1}).
\begin{equation}
{F_{{\mu}}}(x,y) = \sum\nolimits_{i \in {{\cal V}_1}} {{f_{1i{\mu}}}(x,y)}, \label{2-1-1}
\end{equation}
where ${f_{1i{\mu }}}(x,y) = \frac{1}{{{\kappa _1}{\kappa _2}}}\int\nolimits_{\cal U} \int\nolimits_{\cal W} {f_{1i}}(x + {\mu }\xi ,y + {\mu }\eta ){e^{ - \frac{1}{2}{{\left\| \xi  \right\|}^2}}}$ $\cdot {e^{ - \frac{{1}}{2}{\left\| \eta  \right\|}^2}}d\eta d \xi $ is the two-variable Gaussian approximation of ${f_{1i}}(x,y)$, ${\mu } \ge 0$ is the smooth parameters, ${\kappa _1}=\int_{\cal U} {{e^{ - \frac{1}{2}{{\left\| \xi  \right\|}^2}}}d} \xi $, ${\kappa _2}=\int_{\cal W} {{e^{ - \frac{1}{2}{{\left\| \eta  \right\|}^2}}}d} \eta $, and ${\cal U} \subseteq {{\mathbb R}^{{n_1}}}$ and ${\cal W} \subseteq {{\mathbb R}^{{n_2}}}$ are two sets that are symmetric about $0_{n_1}$ and $0_{n_2}$, respectively. In addition,
\begin{equation*}
{ \tilde F_{{\mu}}}(x,y) = - \sum\nolimits_{j \in {{\cal V}_2}} {{f_{2j{\mu}}}(x,y)},
\end{equation*}
where ${f_{2j{\mu }}}(x,y) = \frac{1}{{{\kappa _1}{\kappa _2}}}\int\nolimits_{\cal U} \int\nolimits_{\cal W} {f_{2j}}(x+{\mu }\xi ,y + {\mu }\eta ) {e^{ - \frac{1 }{2}{{\left\| \xi  \right\|}^2}}}$ $\cdot {e^{ - \frac{1}{2}{{\left\| \eta  \right\|}^2}}}d\eta d \xi $ is the  two-variable Gaussian approximation of  ${f_{2j}}(x,y)$. It follows from (\ref{1-1-1-1}) that
\begin{equation}
{F_{{\mu}}}(x,y) = \sum\nolimits_{j \in {{\cal V}_2}} {{f_{2j{\mu}}}(x,y)}. \label{2-2-1}
\end{equation}
Denote ${\cal Z}_\mu ^*={\cal X}_\mu ^* \times {\cal Y}_\mu ^*$ as the nonempty set composed of the saddle point  $(x_\mu ^*,y_\mu ^*)$ of $F_\mu (x,y)$ on ${\cal X} \times {\cal Y}$.

Note that in (\ref{2-1-1}) and (\ref{2-2-1}), the differentiability of ${f_{ \iota s{\mu }}}(x,y)$, $ \iota =1,2$, $ s \in {\cal V}_\iota $, with respect to $x$ and $y$ depends on $\mu $. For $\mu >0$, it is not hard to  see that ${f_{ \iota s{\mu }}}(x,y)$ is differentiable with respect to both $x$ and $y$, i.e., the gradients ${\nabla _x}{f_{\iota s\mu }}( \cdot ,y)$ and ${\nabla _y}{f_{\iota s \mu }}(x, \cdot )$ exist, where $ \iota =1,2$, $ s \in {\cal V}_\iota $.

%Before presenting the distributed stochastic subgradient-free algorithm, we first make the following assumptions.

\begin{assumption}
\label{ass1}
$f_{{\iota}s}(x,y)$, $\iota=1,2$, $s \in {\cal V}_{\iota}$, is convex-concave and locally Lipschitz on a set pair containing ${\cal X} \times {\cal Y}$. The subgradients of $f_{{\iota}s}$ are bounded on the set pair containing ${\cal X} \times {\cal Y}$, i.e., there exist positive constants $D_1$ and $D_2$ such that $\left\| {{\partial _x}{f_{{\iota s}}}(x,y)} \right\| \le D_1$ and $\left\| {{\partial _y}{f_{\iota s}}(x,y)} \right\| \le D_2$ for all $\iota \in \{1,2 \}$ and $s \in {\cal V}_{\iota}$. For each $\iota \in \{1,2\}$, at least one of $f_{\iota s}(x,y)$, $s \in {\cal V}_{\iota }$, is strictly convex-concave.
\end{assumption}

%Let $Z^*= X^* \times Y^*\subseteq {\cal X} \times {\cal Y}$ and $Z_{{\mu }}^*= X_{{\mu}}^* \times Y_{{\mu }}^*\subseteq {\cal X} \times {\cal Y}$ as the sets of all saddle points of $F(x,y)$  on ${\cal X} \times {\cal Y}$ and $F_{{\mu }}(x,y)$ on ${\cal X} \times {\cal Y}$, respectively.

%\begin{assumption}
%\label{ass1}
%$F(x,y)$ has at least a saddle point $(x^*,y^*)$ on ${\cal X} \times {\cal Y}$, i.e., $Z^* \ne \emptyset $. $F_{{{\mu }}}(x,y)$ has at least a saddle point  $(x_{{{\mu }}}^*,y_{{{\mu}}}^*)$ on ${\cal X} \times {\cal Y}$, i.e., $Z_{{{\mu }}}^* \ne \emptyset $.
%\end{assumption}

\begin{assumption}
\label{ass2}
The directed graphs ${\cal G}_1$ and ${\cal G}_2$ are strongly connected and unbalanced.
\end{assumption}

\subsection{The proposed algorithm}

%\subsection{Algorithm description}

To achieve the Nash equilibrium seeking of the considered two-network zero-sum games, we propose the distributed stochastic subgradient-free algorithm. In the proposed algorithm, each agent $i \in {{\cal V}_1}$ maintains two variables $x^{k}_i \in {\mathbb R}^{n_1}$ and $\upsilon ^{k }_i  \in {\mathbb R}^{m_1}$, and each agent $j \in {{\cal V}_2}$ maintains two variables $y^{k}_j \in {\mathbb R}^{n_2}$ and $\omega ^{k }_j  \in {\mathbb R}^{m_2}$. To be more specific, $x^{k}_i$ and $y^{k}_j$ are respectively the states of agents $i $ and $j$, and $\upsilon ^{k }_i$ and $\omega^{k }_j$ are the auxiliary variables. At $k$-th step, where $k=0,1,\ldots$, the variables of agents $i$ and $j$ are updated according to
\begin{align}
x^{k + 1}_i &= {P_{\cal X}}\Big(\sum\nolimits_{s \in {{\cal V}_1}} {{a_{1is}}x_s^k}  - \frac{{{\alpha _k}}}{{\upsilon ^k_{ii}}}g_{1i{\mu}}(x^k_i,\Pi _{1i}^{{k}})\Big),  \quad \notag \\
\upsilon^{k + 1}_i &= \sum\nolimits_{s \in {{\cal V}_1}} {{a_{1is}}}\upsilon^k_s,  \notag \\
y^{k + 1}_j &= {P_{\cal Y}}\Big(\sum\nolimits_{s \in {{\cal V}_2}} {{a_{2js}}y_s^k}  + \frac{{{\alpha _k}}}{{\omega ^k_{jj}}}{g_{2j\mu}}(\Pi _{2j}^{{k}},y^k_j)\Big),  \quad \notag \\
\omega^{k + 1}_j &= \sum\nolimits_{s \in {{\cal V}_2}} {{a_{2js}}}\omega^k_s,  \quad \label{2}
\end{align}
where ${{a_{1is}}}$ and ${{a_{2js}}}$ are the entries of the adjacency matrices $A_1$ and $A_2$. In addition, $\upsilon ^k_i = {(\upsilon ^k_{i1}, \ldots ,\upsilon ^k_{i{m_1}})^T}$, $i \in {{\cal V}_1}$, $\omega ^k_j = {(\omega^k_{j1}, \ldots ,\omega ^k_{j{m_2}})^T}$, $j \in {{\cal V}_2}$. $g_{1i{\mu}}(x^k_i,\Pi _{1i}^{{k}})$ and $g_{2j{\mu }}(\Pi _{2j}^{{k}},y^k_j)$ are the stochastic subgradient-free two-variable oracles which are defined as:
\begin{eqnarray*}
& &{g_{1i\mu }}(x^k_i,\Pi _{1i}^{{k}}) \notag \\
&=& \frac{{f_{1i}}(x^k_i + {\mu }\xi_{1i}^k ,\Pi _{1i}^{{k}} + {\mu }\eta_{1i}^k ) - {f_{1i}}(x^k_i,\Pi _{1i}^{{k}})}{{{\mu }}}\xi_{1i}^k,
\end{eqnarray*}
\begin{eqnarray*}
& &{g_{2j\mu}}(\Pi _{2j}^{{k}},y^k_j)  \notag \\
&=&\frac{{f_{2j}}(\Pi _{2j}^{{k}} + {\mu}\xi_{2j}^k ,y^k_j + {\mu}\eta_{2j}^k )- {f_{2j}}(\Pi _{2j}^{{k}},y^k_j)}{{{\mu }}}\eta_{2j}^k,
\end{eqnarray*}
%\begin{equation*}
%{g_{1i\mu }}(x^k_i,\Pi _{1i}^{{k}})= \frac{{f_{1i}}(x^k_i + {\mu }\xi_{1i}^k ,\Pi _{1i}^{{k}} + {\mu }\eta_{1i}^k ) - {f_{1i}}(x^k_i,\Pi _{1i}^{{k}})}{{{\mu }}}\xi_{1i}^k,
%\end{equation*}
%\begin{equation*}
%{g_{2j\mu}}(\Pi _{2j}^{{k}},y^k_j) =\frac{{f_{2j}}(\Pi _{2j}^{{k}} + {\mu}\xi_{2j}^k ,y^k_j + {\mu}\eta_{2j}^k )- {f_{2j}}(\Pi _{2j}^{{k}},y^k_j)}{{{\mu }}}\eta_{2j}^k,
%\end{equation*}
where $\mu >0$, $\xi_{1i}^k$, $\xi_{2j}^k \in {\mathbb R}^{n_1}$ and $\eta_{1i}^k$, $\eta_{2j}^k \in {\mathbb R}^{n_2}$ are random variables that are generated from the standard Gaussian distribution over ${\cal U}$ and  $ {\cal W}$, respectively, $\Pi _{1i}^{{k}}$, $i \in {\cal V}_1$, and  $\Pi _{2j}^{{k}}$, $j \in {{\cal V}_2}$, will be determined later. ${\alpha _k}$ is a step size satisfying ${\alpha _k} > 0$, $\sum\nolimits_{k=0}^\infty  {{\alpha _k}}  = \infty$ and $\sum\nolimits_{k=0}^\infty  {\alpha _k^2}  < \infty$.
The agents $i  \in {\cal V}_1$ and $j \in {\cal V}_2$ initialize the variables as $x^0_i \in {\cal X}$, $y^0_j \in {\cal Y}$, $\upsilon^{0}_i ={\rm I}_i$, and $\omega^{0}_j ={\rm I}_j$.

Each agent $i \in {\cal V}_1$ and $j \in {\cal V}_2$ perform the first and the third equations of (\ref{2}) to compute the two vectors involved in the cost functions, i.e., $x$ and $y$, respectively. $\Pi _{1i}^{{k}}$, $i \in {\cal V}_1$, is the estimation of the vector $y$ by agent $i$ at step $k$ and $\Pi _{2j}^{{k}}$, $j \in {{\cal V}_2}$, is the estimation of the vector $x$ by agent $j$  at step $k$. According to \cite{Lou16,yang18}, we can set $\Pi _{1i}^{{k}}= \frac{1}{ \sum\nolimits_{s \in {{\cal V}_2}} {{a_{1is}} }}\sum\nolimits_{s \in {{\cal V}_2}} {{a_{1is}}y_s^k} $ and $\Pi _{2j}^{{k}}=\frac{1}{ \sum\nolimits_{s \in {{\cal V}_1}} {{a_{2js}}} } \sum\nolimits_{s \in {{\cal V}_1}} {{a_{2js}}x_s^k} $, where ${a_{1is}} >0$ if ${s \in {{\cal V}_2}}$ is a neighbor of $i$ and ${a_{1is}} =0$ otherwise, in addition, ${a_{2js}}$ is similarly defined for ${s \in {{\cal V}_1}}$. The second and the fourth equations  of (\ref{2}) provide the estimations of the left Perron eigenvectors of $A_1$ and $A_2$ by agents $i  \in {\cal V}_1$ and $j  \in {\cal V}_2$, respectively. In addition, ${{\upsilon ^k_{ii}}}$ and ${\omega ^k_{jj}}$ are used to scale the stochastic subgradient-free two-variable oracles such that the imbalance issue is properly handled since  $A_1$ and $A_2$ are only row-stochastic. According to Perron-Frobenius theorem \cite{Berma94}, $A_\iota $, $\iota =1,2$, has a positive left eigenvector ${\rho _\iota } = {({\rho _{\iota 1}}, \ldots ,{\rho _{\iota {m_\iota }}})^T}$ corresponding to the eigenvalue $1$ and satisfying $\sum\nolimits_{s \in {{\cal V}_\iota }} {{\rho _{\iota s}}}  = 1$.

%Under Assumption {\ref{ass2}} and the fact that the diagonal entries of the row-stochastic matrix $A_\iota $, $\iota =1,2$, are positive, we know from the Perron-Frobenius theorem \cite{Berma94} that $A_\iota $, $\iota =1,2$, has a positive left eigenvector ${\rho _\iota } = {({\rho _{\iota 1}}, \ldots ,{\rho _{\iota {m_\iota }}})^T}$ corresponding to the eigenvalue $1$ and satisfying $\sum\nolimits_{s \in {{\cal V}_\iota }} {{\rho _{\iota s}}}  = 1$.

%Moreover, we note that the neighbors' out-degree of each agent is not used in (\ref{2}).

Let ${{\bar x}^k}$ and ${{\bar y}^k}$ be defined  as follows:
\begin{equation}
{{\bar x}^k} = \sum\nolimits_{i \in {{\cal V}_1}} {{\rho _{1i}}x_i^k}, \ \ {{\bar y}^k} = \sum\nolimits_{j \in {{\cal V}_2}} {{\rho _{2j}}y_j^k}. \label{4.3}
\end{equation}

\begin{assumption}
\label{ass3}
For any fixed $\iota \in \{1,2\}$ and $s \in {\cal V}_\iota $, both the sequences $\{\xi _{\iota s}^k\}$ and $\{\eta _{\iota s}^k\}$ are independent and identically distributed. For any $i_1, i_2, i_3, i_4  \in {\cal V}_1$, $j_1, j_2, j_3, j_4  \in {\cal V}_2$, $i_1 \ne  i_2$, $i_3 \ne  i_4$, $j_1 \ne  j_2$ and $j_3 \ne  j_4$, the  sequences  $\{\xi _{1i_1}^k\}$, $\{\xi _{1i_2}^k\}$, $\{\eta _{1i_3}^k\}$, $\{\eta _{1i_4}^k\}$, $\{\xi _{2j_1}^k\}$, $\{\xi _{2j_2}^k\}$, $\{\eta _{2j_3}^k\}$ and $\{\eta _{2j_4}^k\}$ are mutually independent.
\end{assumption}

The following lemma presents the fundamental properties of the  two-variable Gaussian approximation functions and the stochastic subgradient-free two-variable oracles. The proof of this lemma is given in Appendix A. Let ${\cal F}_k$ be the $\sigma $-field generated by the entire history of the random variables  in algorithm (\ref{2}) from step $0$ to $k-1$. Specifically, ${{\cal F}_0} = \{ x_i^0,y_j^0,i \in {{\cal V}_1},j \in {{\cal V}_2}\}$, ${{\cal F}_k} = \{ x_i^0,y_j^0;{\xi_{1i} ^r},{\eta_{1i} ^r},{\xi_{2j} ^r},{\eta_{2j} ^r};i \in {{\cal V}_1},j \in {{\cal V}_2},0 \le r \le k - 1\}$ for $k \ge 1$.

%The proof of this lemma is given in Appendix A.

\begin{lemmax}
\label{lem1}
Suppose that Assumptions {\ref{ass1}} and  {\ref{ass3}} hold. Let $\tilde {\cal X}^k = \{ {x_\mu^*},x_i^k,{\bar x^k},\Pi _{2j}^k | i \in {\cal V}_1,  j \in {\cal V}_2  \}$, $\tilde {\cal Y}^k = \{ {y_\mu^*},$ $y_j^k,{\bar y^k},\Pi _{1i}^k|  i \in {\cal V}_1, j \in {\cal V}_2\}$. Then,
\begin{enumerate}
\item[\textrm{(i)}] $f_{\iota s{\mu}}(x,y)$ is convex-concave on ${\cal X} \times {\cal Y}$ for any $\iota  \in \{1,2\}$ and $s \in {\cal V}_\iota$. In addition, for each $\iota \in \{1,2\}$, at least one of $f_{\iota s{\mu}}(x,y)$ is strictly convex-concave on ${\cal X} \times {\cal Y}$, and
$f_{\iota s}(x,y) - c \le f_{\iota s\mu}(x,y) \le f_{\iota s}(x,y) + d$ for any $ \iota  \in \{1,2\}$ and $s \in {\cal V}_\iota,$
where $c=\mu D_2 \sqrt{n_2}$, $d=\mu( D_1 \sqrt{n_1} +D_2 \sqrt{n_2})$, $D_1$ and $D_2$ are given in Assumption {\ref{ass1}}, and $n_1$ and $n_2$ are the dimensions of the vectors $x$ and $y$, respectively;
\item[\textrm{(ii)}] for any $i \in {\cal V}_1$, $j \in {\cal V}_2$ and $ k \ge 0$, it holds that ${\nabla _{x_i^k}}{f_{1i\mu }}(x_i^k,\Pi _{1i}^k) = \Xi [{g_{1i\mu }}(x_i^k,\Pi _{1i}^k)|{{\cal F}_k}]$, ${\nabla _{y_j^k}}{f_{2j\mu }}(\Pi _{2j}^k,y_j^k) = \Xi [{g_{2j\mu }}(\Pi _{2j}^k,y_j^k)$ $\big|{  {\cal F}_k} ]$;
\item[\textrm{(iii)}]  for any $\iota \in \{1,2\}$, $s  \in {\cal V}_\iota$, $(\tilde x^k,\tilde y^k) \in \tilde {\cal X}^k \times \tilde {\cal Y}^k $ and $ k \ge 0$, it holds that $\left\| {{\nabla _{{{\tilde x}^k}}}{f_{\iota s \mu }}({{\tilde x}^k},{{\tilde y}^k})} \right\| \le {M_1}$, $\left\| {{\nabla _{{{\tilde y}^k}}}{f_{\iota s \mu }}({{\tilde x}^k},{{\tilde y}^k})} \right\| \le {M_2}$, where ${M_\iota } = [D_\iota ^2{(4 + {n_\iota })^2} + D_{3 - \iota }^2{n_1}{n_2} +2 {D_1}{D_2}{(3 + {n_\iota })^{\frac{3}{2}}}({n_{3 - \iota }})^{\frac{1}{2}} ]^{\frac{1}{2}}$, $\iota  = 1,2$;
\item[\textrm{(iv)}] for any $i \in {\cal V}_1$, $j \in {\cal V}_2$ and $ k \ge 0$, it holds that $\Xi[{\left\| {g_{1i{\mu }}(x_i^k,\Pi _{1i}^k)} \right\|^2}|{{\cal F}_k}] \le M_1^2$, $\Xi [{ {\left\| {g_{1i\mu}(x_i^k,\Pi _{1i}^k)} \right\|^2} }] \le {M^2_1}$, $\Xi[{\left\| {g_{2j{\mu }}(\Pi _{2j}^k,y_j^k)} \right\|^2}|{{\cal F}_k}] \le M_2^2$, and $\Xi [ \| g_{2j\mu}(\Pi _{2j}^k,$ $y_j^k) \|^2 ] \le {M^2_2}$, where ${M_\iota }$, $\iota  = 1,2$, is given in (iii).
\end{enumerate}
\end{lemmax}
\medskip

\section{Main results}

In this section, we first show  the almost sure consensus  of the agents' states for each subnetworks, i.e., ${\cal V}_1$ and ${\cal V}_2$. We then establish the almost sure convergence to the Nash equilibrium of the agents' states of the overall network.

To facilitate the subsequent analysis, we rewrite the first and the third equations in (\ref{2}) as
\begin{eqnarray}
x_i^{k + 1} &=&  \sum\nolimits_{s \in {{\cal V}_1}} {{a_{1is}}x_s^k + \varepsilon _{1i}^k}, \notag \\
y_j^{k + 1} &=& \sum\nolimits_{s \in {{\cal V}_2}} {{a_{2js}}y_s^k}  + \varepsilon _{2j}^k, \label{1-1}
\end{eqnarray}
where $\varepsilon _{1i}^k = {P_{\cal X}}(\sum\nolimits_{s \in {{\cal V}_1}} {{a_{1is}}x_s^k}  - \frac{{{\alpha _k}}}{{\upsilon _{ii}^k}}{g_{1i\mu }}(x_i^k,\Pi _{1i}^k)) - \sum\nolimits_{s \in {{\cal V}_1}}{{a_{1is}}x_s^k}$, $\varepsilon _{2j}^k = {P_{\cal Y}}(\sum\nolimits_{s \in {{\cal V}_2}} {{a_{2js}}y_s^k}  + \frac{{{\alpha _k}}}{{\omega _{jj}^k}}{g_{2j\mu }}(\Pi _{2j}^k,y_j^k))$ $ - \sum\nolimits_{s \in {{\cal V}_2}} {{a_{2js}}y_s^k} $.

%where $\varepsilon _{1i}^k = {P_{\cal X}}(\sum\limits_{s \in {{\cal V}_1}} {{a_{1is}}x_s^k}  - \frac{{{\alpha _k}}}{{\upsilon _{ii}^k}}{g_{1i\mu }}(x_i^k,\Pi _{1i}^k)) - \sum\limits_{s \in {{\cal V}_1}}{{a_{1is}}x_s^k}$, $\varepsilon _{2j}^k = {P_{\cal Y}}(\sum\limits_{s \in {{\cal V}_2}} {{a_{2js}}y_s^k}  + \frac{{{\alpha _k}}}{{\omega _{jj}^k}}{g_{2j\mu }}(\Pi _{2j}^k,y_j^k)) - \sum\limits_{s \in {{\cal V}_2}} {{a_{2js}}y_s^k} $.

\subsection{Almost sure consensus}

%The following Lemma shows the accumulation error between the agents' states and the weighted states accompanied by the step size.

%We first present a preliminary result regarding $\varepsilon _{1i}^k $ and $\varepsilon _{2j}^k$, $i  \in {\cal V}_1$, $j \in {\cal V}_2$.

We first present a preliminary result and then establish the almost sure consensus of the agents' states for each subnetworks. The proofs of Lemma {\ref{lem4}} and Theorem {\ref{thm1}} are given in Appendices B and C, respectively.

\begin{lemmax}
\label{lem4}
For all $i  \in {\cal V}_1$,  $j \in {\cal V}_2$ and any  $k \ge 0$, it holds that $\left\| {\varepsilon _{1i}^k} \right\| \le \frac{{{\alpha _k}}}{{\upsilon _{ii}^k}}\left\| {g_{1i \mu }(x_i^k,\Pi _{1i}^k)} \right\|$, and $\left\| {\varepsilon _{2j}^k} \right\| \le \frac{{{\alpha _k}}}{{\omega _{jj}^k}}\| g_{2j \mu }(\Pi _{2j}^k, y_j^k) \|$.
\end{lemmax}
\medskip

\begin{theoremx}
\label{thm1}
Suppose that Assumptions {\ref{ass1}}-{\ref{ass3}} hold. Let ${\{ x_i^k\} _{k \ge 0}}$ and ${\{ y_j^k\} _{k \ge 0}}$ be generated by (\ref{2}). Then,  all agents in ${\cal V}_1$ achieve consensus almost surely, while all agents in ${\cal V}_2$ achieve consensus almost surely. Specifically, for each $i  \in {\cal V}_1$, $\mathop {\lim }\nolimits_{k \to \infty } \| x_i^k - {\bar x^k} \|= 0$  with probability one, and  for each $ j \in {\cal V}_2$, $\mathop {\lim }\nolimits_{k \to \infty } \| y_j^k - {\bar y^k} \|  = 0$ with probability one, where $\bar x^k$ and $\bar y^k$ are given in  (\ref{4.3}).
\end{theoremx}
\medskip

\subsection{Almost sure convergence to the Nash equilibrium}

The following lemma presents the relation between $\Xi [ $ $(\|{\bar x}^{k + 1} { {{} - x_{{\mu }}^*} \|^2}+{\| {{{\bar y}^{k + 1}} - y_{{\mu }}^*} \|^2})  \big|{{\cal F}_k}]$ and $\| {{{\bar x}^k} - x_{{\mu }}^*} \|^2+ \| {\bar y}^k $ $- y_{{\mu }}^* \|^2$, which will be used to show the almost sure convergence.  The proof of this lemma is given in Appendix D.

\begin{lemmax}
\label{lem7}
Suppose that Assumptions {\ref{ass1}}-{\ref{ass3}} hold. Then, for any $(x^*_{\mu},y^*_{\mu}) \in {\cal Z}^*_{\mu}$ and $k \ge 0$, it holds that
\begin{eqnarray}
& &\Xi [ ({\left\| {{{\bar x}^{k + 1}} - x_{{\mu }}^*} \right\|^2}+{\left\| {{{\bar y}^{k + 1}} - y_{{\mu }}^*} \right\|^2})  |{{\cal F}_k}]  \notag \\
&\le & (1 + \alpha _k^2)({\left\| {{{\bar x}^k} - x_{{\mu }}^*} \right\|^2+\left\| {{{\bar y}^k} - y_{{\mu }}^*} \right\|^2})+ p^k-q^k,  \qquad \label{25}
\end{eqnarray}
where $p^k$ is a nonnegative random variable satisfying $\sum\nolimits_{k = 0}^\infty  {\Xi [{p^k}]} < \infty$ and its expression is given in  (\ref{36}), $q^k=2{\alpha _k}[ -F_{{\mu }}(x_{{\mu }}^*,{{\bar y}^k}) + F_{{\mu }}(x_{{\mu }}^*,y_{{\mu }}^*) + F_{{\mu }}({{\bar x}^k},$ $y_{{\mu}}^*)- F_{{\mu}}(x_{{\mu }}^*,y_{{\mu }}^*)]$.
\end{lemmax}
\medskip

To this end, we have established the convergence result that the agents' states  converge almost surely to the neighborhood of the Nash equilibrium, and also specify the bounds of the differences between the limits of the cost values of all the agents and the cost value at the Nash equilibrium.

%?????????like Hu Guoqiang, in fact, he derive the almost sure convergence about the state but he give the convergence of the $f$ in expectation, so I can do like this to make this theorem is in expectation

\begin{theoremx}
\label{thm2}
Suppose that Assumptions {\ref{ass1}}-{\ref{ass3}} hold. Let $(x^*,y^*)$ be the Nash equilibrium, and ${\{ x_i^k\} _{k \ge 0}}$ and ${\{ y_j^k\} _{k \ge 0}}$ be generated by (\ref{2}). Then, for each $i \in {\cal V}_1$, $\mathop {\lim }\nolimits_{k \to \infty } F({ x_i^k},{y^*}) -F({x^*},{y^*}) \le \vartheta$  with probability one, and for each $j \in {\cal V}_2$, $\mathop {\lim }\nolimits_{k \to \infty } {F}(x^*,{y_j^k}) -{F}(x^*,y^*) \ge -  \vartheta$ with probability one, where $\vartheta =\mu \min \{m_1,m_2\}( D_1 \sqrt{n_1} +2D_2 \sqrt{n_2})$.
\end{theoremx}
\noindent \textbf{Proof:}\,\, Consider the equation (\ref{25}) given in Lemma {\ref{lem7}}. Since $p^k \ge 0$ for all $k \ge 0$, $\{\sum\nolimits_{k = 0}^\ell  {{p^k}}\}_{\ell  \ge 0}$ is a nondecreasing sequence.  Then, it follows that  $\Xi [\sum\nolimits_{k = 0}^\infty  {{p^k}} ] = \Xi [\mathop {\lim }\limits_{\ell  \to \infty } \sum\nolimits_{k = 0}^\ell  {{p^k}} ] = \mathop {\lim }\limits_{\ell  \to \infty } \Xi [\sum\nolimits_{k = 0}^\ell  {{p^k}} ] =\mathop {\lim }\limits_{\ell  \to \infty } \sum\nolimits_{k = 0}^\ell  \Xi  [{p^k}] = \sum\nolimits_{k = 0}^\infty  \Xi  [{p^k}]$, where the monotone convergence theorem has been used in the second equality. Then, using the fact that $\sum\nolimits_{k = 0}^\infty  {\Xi [{p^k}]} < \infty$, we know that $\Xi [\sum\nolimits_{k = 0}^\infty  {{p^k}} ] < \infty $. This implies that $ \sum\nolimits_{k = 0}^{\infty} p^k < \infty $ a.s.. Moreover, we know that $q^k \ge 0$ for any $k \ge 0$ since $ F_{\mu}(x_{\mu}^*,{{\bar y}^k}) \le F_{\mu}(x_{\mu}^*,y_{\mu}^*)$ and $ F_{\mu}({{\bar x}^k},y_{\mu}^*) \ge F_{\mu}(x_{\mu}^*,y_{\mu}^*) $ for any $k \ge 0$, and $\sum\nolimits_{k = 0}^\infty  {\alpha _k^2}  < \infty$. Therefore, according to the supermartingale convergence lemma due to Robbins-Siegmund (see Lemma 11 in \cite{polyak87}), we can obtain that ${\left\| {{{\bar x}^{k}} - x_{\mu}^*} \right\|^2}+{\left\| {{{\bar y}^{k }} - y_{\mu}^*} \right\|^2}$ converges to a certain nonnegative variable a.s. for any $(x_{\mu}^*,y_{\mu}^*) \in {\cal Z}_{\mu}^*$, and $\sum\nolimits_{k = 0}^\infty  {{q^k}}  < \infty $ a.s.. The first result  implies that both ${{\bar x}^{k}} $ and ${{\bar y}^{k}}$ are bounded a.s., and hence the sequences $\{{{\bar x}^{k}} \}$ and $\{{{\bar y}^{k}}\}$ have convergent subsequences a.s.. The second result  and the fact $\sum\nolimits_{k = 0}^\infty  {{\alpha _k} = \infty } $ imply that $\mathop {\lim \inf }\limits_{k \to \infty } [ -F_{\mu}(x_{\mu}^*,{{\bar y}^k}) + F_{\mu}(x_{\mu}^*,y_{\mu}^*) + F_{\mu}({{\bar x}^k},y_{\mu}^*)- F_{\mu}(x_{\mu}^*,y_{\mu}^*)]  =0$  with probability one. It then follows that $\mathop {\lim \inf }\nolimits_{k \to \infty }  {F_{\mu}}(x_{\mu}^*,{{\bar y}^{{k}}})= {F_{\mu}}(x_{\mu}^*,y_{\mu}^*)$ and $\mathop {\lim \inf }\nolimits_{k \to \infty }{F_{\mu}}({{\bar x}^{{k}}},y_{\mu}^*) = {F_{\mu}}(x_{\mu}^*,y_{\mu}^*)$ with probability one. Note that ${F_{\mu}}$ is continuous. Then, we know that there exist convergent subsequences $\{{{\bar x}^{k_1}}\}$ and $\{{{\bar y}^{k_1}}\}$ such that $\mathop {\lim }\nolimits_{k_1 \to \infty }  {F_{\mu}}(x_{\mu}^*,{{\bar y}^{{k_1}}})={F_\mu }(x_\mu ^*,{\lim _{{k_1} \to \infty }}{\bar y^{{k_1}}})= {F_{\mu}}(x_{\mu}^*,y_{\mu}^*)$ and $\mathop {\lim }\nolimits_{k_1 \to \infty }{F_{\mu}}({{\bar x}^{{k_1}}},y_{\mu}^*) ={F_\mu }({\lim _{{k_1} \to \infty }}{\bar x^{{k_1}}},y_\mu ^*)= {F_{\mu}}(x_{\mu}^*,y_{\mu}^*)$ with probability one. Furthermore, we know that $\mathop {\lim }\nolimits_{{k_1} \to \infty } {\bar x^{{k_1}}} = \tilde x_{\mu}$ and $\mathop {\lim }\nolimits_{{k_1} \to \infty } {\bar y^{{k_1}}} = \tilde y_{\mu}$  for some pair $(\tilde x_{\mu},\tilde y_{\mu}) \in {\cal Z}_{\mu}^* $ with probability one. According to Lemma {\ref{lem1}} (i), ${F_{\mu}}$ is strictly convex-concave, then it has a unique saddle point on ${\cal X} \times {\cal Y}$, i.e., ${\cal Z}^*_{\mu}=\{(x^*_{\mu},y^*_{\mu}) \}$. Therefore, $\mathop {\lim }\nolimits_{{k_1} \to \infty } {\bar x^{{k_1}}} = x_{\mu}^*$ and $\mathop {\lim }\nolimits_{{k_1} \to \infty } {\bar y^{{k_1}}} = y_{\mu}^*$  with probability one. This together with  the fact that ${\| {{{\bar x}^{k}} - x_{\mu}^*} \|^2}+{\| {{{\bar y}^{k }} - y_{\mu}^*} \|^2}$ is convergent a.s. implies that ${\| {{{\bar x}^{k}} - x_{\mu}^*} \|^2}+{\| {{{\bar y}^{k }} - y_{\mu}^*} \|^2}$ converges to zero a.s.. This further indicates that ${\lim _{k \to \infty }}\|{{\bar x}^k} - x_\mu ^*\|=0$ (or ${\lim _{k \to \infty }}\|{{x}_i^k} - x_\mu ^*\|=0$, $ i \in {\cal V}_1$) with probability one and ${\lim _{k \to \infty }}\|{{\bar y}^k} - y_\mu ^*\|=0$ (or ${\lim _{k \to \infty }}\|{{y}_j^k} - y_\mu ^*\|=0$, $ j \in {\cal V}_2$) with probability one.

Let $({x^*},{y^*}) \in {\cal Z}^*$. Clearly, $({x^*},{y^*})$ is exactly the Nash equilibrium of the two-network zero-sum game. By virtue of the continuity of $F_{{\mu }}$, one has $\mathop {\lim }\nolimits_{k \to \infty } {F_{{\mu }}}({{\bar x}^k},{y^*}) = {F_{{\mu }}}({x_{{\mu}}^*},y^*)$ with probability one. Denote $\vartheta_1 =\mu \min \{m_1,m_2\}D_2 \sqrt{n_2} $, $\vartheta_2 =\mu \min \{m_1,m_2\}( D_1 \sqrt{n_1} +D_2 \sqrt{n_2})$. Based on  Lemma {\ref{lem1}} (i), we know that $F({{\bar x}^k},{y^*}) \le {F_{{\mu }}}({{\bar x}^k},{y^*}) + \vartheta_1$  for any $k \ge 0$, and ${F_{{\mu }}}(x_{{\mu}}^*,{y^*}) \le F_{{\mu }}(x_{{\mu}}^*,{y_{{\mu }}^*}) \le F_{{\mu }}(x^*,{y_{{\mu }}^*})\le F(x^*,{y_{{\mu }}^*}) + \vartheta_2 \le F(x^*,{y^*}) + \vartheta_2$. Therefore, $\mathop {\lim }\nolimits_{k \to \infty } F({{\bar x}^k},{y^*}) \le F({x^*},{y^*}) + \vartheta$ with probability one, where $ \vartheta= \vartheta_1 + \vartheta_2$.  Furthermore, for each $i \in {\cal V}_1$, one has $\mathop {\lim }\limits_{k \to \infty } F(x_i^k,{y^*}) \le \mathop {\lim }\limits_{k \to \infty } F({{\bar x}^k},{y^*}) + | \mathop {\lim }\limits_{k \to \infty } ({F(x_i^k,{y^*}) - F({{\bar x}^k},{y^*})} )|  \le F({x^*},{y^*}) + \vartheta$ with probability one, where one has used Theorem {\ref{thm1}} in the second inequality. On the other hand, it follows from the continuity of $F_{{\mu }}$ that $\mathop {\lim }\nolimits_{k \to \infty } {F_{{\mu }}}({x^*},{{\bar y}^k}) = {F_{{\mu }}}({x^*},{y_{{\mu}}^*})$ with probability one. Based on Lemma {\ref{lem1}} (i), we know that $F(x^*,{{\bar y}^k}) \ge {F_{{\mu }}}(x^*,{{\bar y}^k}) -$
$ \vartheta_2$ for any $k \ge 0$, and ${F_{{\mu }}}(x^*,y_{{\mu}}^*) \ge F_{{\mu }}(x_{{\mu }}^*,y_{{\mu }}^*) \ge  F_{{\mu }}(x_{{\mu}}^*,y^*)  $ $ \ge F(x_{{\mu }}^*,y^*) - \vartheta_1  \ge F(x^*,y^*) - \vartheta_1$. Therefore, $\mathop {\lim }\limits_{k \to \infty } {F}(x^*,{\bar y^k})$ $ \ge {F}(x^*,y^*) -  \vartheta$  with probability one. Furthermore, for each  $j \in {\cal V}_2$, $\mathop {\lim }\limits_{k \to \infty } F({x^*},y_j^k) \ge \mathop {\lim }\limits_{k \to \infty } F({x^*},{\bar y^k}) - | \mathop {\lim }\limits_{k \to \infty }(F({x^*},y_j^k)  $ $- F({x^*},{{\bar y}^k}) )| \ge F({x^*},{y^*}) - \vartheta$ with probability one, where one has used  Theorem {\ref{thm1}} in the second inequality.
\medskip

%one has used the fact that $ | \mathop {\lim }\limits_{k \to \infty } (F(x_i^k,{y^*}) - F({{\bar x}^k},{y^*})) |=0$  with probability one in

\begin{remark}
The bound provided in Theorem {\ref{thm1}} is generated due to the use of the subgradient-free two-variable oracles instead of the true subgradients, and it can be sufficiently close to $0$ by choosing the smooth parameter $\mu >0$ to be sufficiently small.
\end{remark}
\medskip

\section{Numerical simulations}

\begin{figure}[t]
\centering
\includegraphics[height=1.7cm ,width=3.8cm]{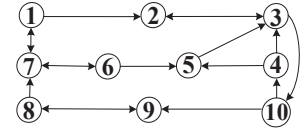}
\captionsetup{margin=20pt,format=hang,justification=justified}
\caption{The graph of the subnetwork}
\label{fig2}
\end{figure}

\begin{figure}[t]
\centering
\includegraphics[height=4cm ,width=5.9cm]{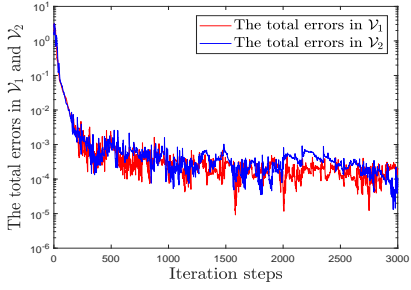}
\caption{The total errors between the agents' states and the Nash equilibrium}
\label{fig3}
\end{figure}

Consider a zero-sum game engaged by two (sub)networks. Each subnetwork has $10$ agents. The sets of the agents are ${\cal V}_1 =\{1, \ldots, 10\}$ and ${\cal V}_2 =\{1, \ldots, 10\}$, respectively. The cost function of the first subnetwork is $F(x,y)=\sum\nolimits_{i \in {\cal V}_1} {{f_{1i}}(x,y)}$, where $f_{11}(x,y)=\frac{1}{{10}}(| x | - | y |) + \cos\frac{y}{2}$ and $f_{1i}(x,y)=\frac{i}{{10}}{x^2} - \frac{i}{5}{y^2} - \frac{1}{{10}}\cos\frac{x}{2}$ for $i \in {\cal V}_1 / \{1\}$. The cost function of the second subnetwork is $\tilde F(x,y)=-\sum\nolimits_{j \in {\cal V}_2} {{f_{2j}}(x,y)}$, where  $f_{21}(x,y)=\frac{1}{{10}}| x | - \frac{1}{{10}}| y | - \frac{9}{{10}}\cos\frac{x}{2} - 9{y^2}$ and $f_{2j}(x,y)=\frac{j}{{10}}{x^2} - \frac{1}{5}{y^2} + \frac{1}{9}\cos\frac{y}{2}$ for $j \in {\cal V}_2 / \{1\}$. Set ${\cal X}={\cal Y}={\cal U}={\cal W}=[-1,1]$. The two subnetworks are equipped with the same strongly connected and unbalanced directed graph shown in Fig. {\ref{fig2}}, where each node has a self-loop which is not drawn in Fig. {\ref{fig2}}. There is an undirected communication link between any $i \in {\cal V}_1$ and $j \in {\cal V}_2$, where $j=i$. It can be shown that $F(x,y)+\tilde F(x,y)=0$ and $F(x,y)$ has a unique saddle point (i.e., the unique Nash equilibrium of the game) $(x^*,y^*)=(0,0)$.  We also set ${\alpha _k} = \frac{0.1}{{k + 1}}$, $\mu=0.001$. Fig. {\ref{fig3}} illustrates  the total errors between the agents' states in the two subnetworks (or ${\cal V}_1$ and ${\cal V}_2$) and the Nash equilibrium $(x^*,y^*)$, i.e., $\sqrt {\sum\nolimits_{i \in {\cal V}_1}{(x_i^k}  - {x^*})^2} $ and $\sqrt {\sum\nolimits_{j \in {\cal V}_2}{(y_j^k}  - {y^*})^2} $. It can be seen from Fig. {\ref{fig3}} that the total errors approach to $10^{-4}$ (i.e., the size of the Nash equilibrium's neighborhood to which all the agents converge) as the iteration steps increase. This is consistent with our theoretical results.

%Nevertheless, we note that the error is quite small which also illustrates the feasibility of the proposed algorithm for achieving the Nash equilibrium of the two-network zero-sum game.

\section{Conclusion}

This paper has developed a distributed stochastic subgradient-free algorithm for achieving the Nash equilibrium in two-network zero-sum games, where the local cost functions are  unknown. The proposed algorithm introduces the subgradient-free two-variable oracles and is applicable to the unbalanced directed graphs. It is shown that by using the proposed algorithm, the agents' states converge  almost surely to a neighborhood of the Nash equilibrium  and the neighborhood can be arbitrary small by taking suitable smooth parameter.

\section*{Appendix A}

For simplicity, we only prove the results for $i \in {\cal V}_1$ and those for $ j  \in {\cal V}_2$ can be similarly obtained.

(i) Based on the definition of $f_{1i{\mu}}$ and the fact that $f_{1i}$ is  (strictly) convex-concave on the set containing ${\cal X} \times {\cal Y}$, it is not hard to show that $f_{1i{\mu}}$ is  (strictly) convex-concave on ${\cal X} \times {\cal Y}$. For any $(x,y) \in {\cal X} \times {\cal Y}$, we denote $\tilde f_i(x,y ) = \frac{1}{{{\kappa _1}{\kappa _2}}}\int\nolimits_{\cal U} \int\nolimits_{\cal W} {\left( {{f_{1i}}(x + {\mu }\xi ,y + {\mu }\eta ) - {f_{1i}}(x,y + {\mu }\eta )} \right)} {e^{ - \frac{1 }{2}{{\left\| \xi  \right\|}^2}}} \cdot $ ${e^{ - \frac{1}{2}{{\left\| \eta  \right\|}^2}}}d\eta  d\xi $, $\tilde {\tilde f}_i(x,y ) = \frac{1}{{{\kappa _1}{\kappa _2}}}\int\nolimits_{\cal U} \int\nolimits_{\cal W} ({f_{1i}}(x,y + {\mu }\eta ) -f_{1i}(x,$ $ y))  {e^{ - \frac{1}{2}{{\left\| \xi  \right\|}^2} }}{e^{ - \frac{1}{2} {{\left\| \eta  \right\|}^2}}}d\eta d\xi $. It is clear that ${f_{1i{\mu }}}(x,y)-{f_{1i}}(x,y)$ $={\tilde f}_i(x,y )+{\tilde {\tilde f}}_i(x,y )$. Let $c$ and $d$ be two constants which are given in Lemma  {\ref{lem1}} (i). Using the convexity of ${f_{1i}}( \cdot ,y + {\mu }\eta )$ at  $x$, one has $\tilde f_i(x,y ) \ge \frac{1}{{{\kappa _1}{\kappa _2}}}\int\nolimits_{\cal W} {(\int\nolimits_{\cal U} {{\mu}{\xi ^T}} } {e^{ - \frac{1}{2}{{{\left\| \xi  \right\|}^2}}}}d\xi ){\partial _x}{f_{1i}}(\cdot,$
$y + {\mu }\eta ){e^{ - \frac{1}{2}{{{\left\| \eta  \right\|}^2}}}}d\eta =0$. Using the concavity of ${{f_{1i}}(x, \cdot )}$ at $y+{\mu }\eta $ and the fact that $\|{{\partial _{y + {\mu}\eta }}{f_{1i}}} (x, \cdot ) \| \le D_2$ given in Assumption {\ref{ass1}}, one has $\tilde {\tilde f}_i (x,y ) \ge  - \frac{1}{{{\kappa _1}{\kappa _2}}}\int\nolimits_{\cal U} \int\nolimits_{\cal W} {\mu}\left\| \eta  \right\|$ $\left\| {{\partial _{y + {\mu }\eta }}{f_{1i}}(x, \cdot )} \right\|  {e^{ - \frac{1}{2}{{\left\| \xi  \right\|}^2} }}{e^{ - \frac{1}{2}{{\left\| \eta  \right\|}^2}}}d\eta d\xi \ge - c$. Therefore, $ f_{1i\mu}(x,y)  \ge f_{1i}(x,y) - c$. By using a symmetric analysis, we can obtain $f_{1i\mu}(x,y) \le f_{1i}(x,y) + d$.

(ii) This part of proof is similar to that of (21) in \cite{Nester17}, which however, concerns the one-variable cost function, i.e., $f_{\mu}(x)$. Suppose that the variables $x \in {\cal X}$, $y  \in {\cal Y}$, $\xi \in {\cal U}$ and $\eta \in {\cal W}$ are independent with respect to each other. Then, all the mathematical manipulations performed for $f_{\mu}(x)$ in the proof of (21) in \cite{Nester17} can be similarly performed for ${f_{1i{\mu }}}(x,y) $. By performing such manipulations and using the facts that ${{f_{1i}}(x ,y )}$ is independent on $(\xi ,\eta)$ and $\int_{\cal U} {{e^{ - \frac{1}{2}{{\left\| \xi  \right\|}^2}}}} \xi d\xi  = 0$, it is not hard to show that ${\nabla _x}{f_{1i\mu }}(x,y) =\frac{1}{{{\kappa _1}{\kappa _2}}}\int_{\cal U} \int_{\cal W}$ $ {
\frac{{{f_{1i}}(x + \mu \xi ,y + \mu \eta )} -{{f_{1i}}(x ,y )}}{\mu }} {e^{ - \frac{1}{2}{{\left\| \xi  \right\|}^2}}}{e^{ - \frac{1}{2} {{\left\| \eta  \right\|}^2}}}\xi d\eta d\xi $, where the detailed calculations are omitted  due to space limitations. Note that for any fixed $k \ge 0$, the variables $x_i^k$, $\Pi _{1i}^k$,  $\xi_{1i}^k$ and $\eta_{1i}^k$ are independent with respect to each other. Moreover, we know that given $i \in {\cal V}_1$, there exists $s \in {\cal V}_2$ such that $a_{1is} >0$. Without loss of generality, we assume $a_{1is} >0$ for each $s \in {\cal V}_2$. Then, based on the facts that $y^k_s \in {\cal Y}$ for all $s \in {\cal V}_2$, ${\cal Y}$ is a convex set, $0<\frac{{{a_{1is}}}}{ \sum\nolimits_{s \in {{\cal V}_2}} {{a_{1is}} }}<1$ and $\frac{1}{ \sum\nolimits_{s \in {{\cal V}_2}} {{a_{1is}} }}\sum\nolimits_{s \in {{\cal V}_2}} {{a_{1is}}}=1$, we  know that $\Pi _{1i}^k \in {\cal Y}$. Thus, it follows that the above obtained expression of ${\nabla _x}{f_{1i\mu }}(x,y)$ holds for $x=x_i^k$, $y=\Pi _{1i}^k$,  $\xi=\xi_{1i}^k$ and $\eta=\eta_{1i}^k$. According to the definition of the conditional expectation ${\Xi[ \cdot|\cdot ]}$ given in Section I and Assumption {\ref{ass3}}, it is not hard to show that ${\nabla _{x_i^k}}{f_{1i\mu }}(x_i^k,\Pi _{1i}^k) = \Xi [{g_{1i\mu }}(x_i^k,\Pi _{1i}^k)|{{\cal F}_k}]$.

(iii) It is not hard to show that $\Pi _{2j}^k \in {\cal X}$, $\bar x^k \in {\cal X}$ and $\bar y^k \in {\cal Y}$. Under Assumption {\ref{ass3}}, it can be shown that the expression of ${\nabla _x}{f_{1i\mu }}(x,y)$ obtained in the proof of (ii) holds for $x=\tilde x^k $, $y=\tilde y^k$, $\xi=\xi_{1i}^k$ and $\eta=\eta_{1i}^k$, where $\tilde x^k  \in \tilde {\cal X}^k$ and $\tilde y^k  \in \tilde {\cal Y}^k $. It then follows that for any fixed $k \ge 0$, $\left\| {{\nabla _{{{\tilde x}^k}}}{f_{1i\mu }}({{\tilde x}^k},{{\tilde y}^k})} \right\| \le \frac{1}{{{\kappa _1}{\kappa _2}}}\int_{\cal U}\int_{\cal W} \sqrt {\tilde g_i^k} {e^{ - \frac{1}{2}{{\left\| {\xi _{1i}^k} \right\|}^2}}}{e^{ - \frac{1}{2}{{\left\| {\eta _{1i}^k} \right\|}^2}}}d\eta _{1i}^k d\xi _{1i}^k=\Xi[\sqrt {\tilde g_i^k}]$, where $\tilde g_i^k = ({D_1}\left\| {\xi _{1i}^k} \right\| + {D_2}\left\| {\eta _{1i}^k} \right\|)^2{\left\| {\xi _{1i}^k} \right\|^2}$. Let ${M_\iota }$, $\iota  = 1,2$, be a constant which is given in Lemma  {\ref{lem1}} (iii). It is shown by the definition of the total expectation ${\Xi[ \cdot ]}$ and Assumption {\ref{ass3}} that ${\Xi[\tilde g_i^k]} \le M^2_1$, where Lemma 1 in \cite{Nester17} has been used. Then, it follows from the relation $\Xi[\sqrt {\tilde g_i^k}] \le \sqrt {\Xi[\tilde g_i^k]} $ that $\left\| {{\nabla _{{{\tilde x}^k}}}{f_{1i\mu }}({{\tilde x}^k},{{\tilde y}^k})} \right\| \le \sqrt {\Xi[{\tilde g_i^k}]} \le  M_1 $. Moreover, using the similar arguments to those employed to obtain ${\nabla _x}{f_{1i\mu }}(x,y)$ in the proof of (ii), we can derive ${\nabla _y}{f_{1i\mu }}(x,y) =\frac{1}{{{\kappa _1}{\kappa _2}}}\int_{\cal U} \int_{\cal W}$ $ {
\frac{{{f_{1i}}(x + \mu \xi ,y + \mu \eta )} -{{f_{1i}}(x ,y )}}{\mu }} {e^{ - \frac{1}{2}{{\left\| \xi  \right\|}^2}}}{e^{ - \frac{1}{2} {{\left\| \eta  \right\|}^2}}}\eta   d\eta d\xi $. It then follows that for any fixed $k \ge 0$, $\left\| {{\nabla _{{{\tilde y}^k}}}{f_{1i\mu }}({{\tilde x}^k},{{\tilde y}^k})} \right\| \le \Xi[\sqrt {{\bar g}_i^k}]$, where ${\bar g}_i^k = ({D_1}\left\| {\xi _{1i}^k} \right\| + {D_2}\left\| {\eta _{1i}^k} \right\|)^2{\left\| {\eta _{1i}^k} \right\|^2}$. It follows from the similar analysis by which we obtain $\left\| {{\nabla _{{{\tilde x}^k}}}{f_{1i\mu }}({{\tilde x}^k},{{\tilde y}^k})} \right\| \le  M_1 $ that $\left\| {{\nabla _{{{\tilde y}^k}}}{f_{1i\mu }}({{\tilde x}^k},{{\tilde y}^k})} \right\| \le {M_2}$.

(iv) Under Assumptions {\ref{ass1}} and {\ref{ass3}}, it is not hard to show that ${\left\| {g_{1i{\mu }}(x_i^k,\Pi _{1i}^k)} \right\|^2} \le \tilde g_i^k $ and $\Xi[\tilde g_i^k|{{\cal F}_k}] =\Xi[\tilde g_i^k]$, where $\tilde g_i^k$ is given in the proof of (iii). Therefore, one has $\Xi[{\left\| {g_{1i{\mu }}(x_i^k,\Pi _{1i}^k)} \right\|^2}|{{\cal F}_k}] \le M_1^2$, $\Xi[{\left\| {g_{1i{\mu }}(x_i^k,\Pi _{1i}^k)} \right\|^2}] \le M_1^2$.
\medskip

\section*{Appendix B}

It is not hard to show that $\sum\nolimits_{s \in {{\cal V}_1}} {{a_{1is}}x_s^k} \in {\cal X}$, $\forall k \ge 0$. From the proof of Proposition 1 in \cite{MaiVS16}, we know that $\upsilon^{k}_{ii}$ is positive for any $k \ge 0$. Then, following Lemma 1 (b) in \cite{Nedic10}, one has $\left\| {\varepsilon _{1i}^k} \right\| \le \frac{{{\alpha _k}}}{{\upsilon _{ii}^k}}\left\| {g_{1i \mu }(x_i^k,\Pi _{1i}^k)} \right\|$ for any  $k \ge 0$. Similarly, we can obtain that $\left\| {\varepsilon _{2j}^k} \right\| \le \frac{{{\alpha _k}}}{{\omega _{jj}^k}}\left\| {g_{2j \mu }(\Pi _{2j}^k,y_j^k)} \right\|$ for any  $k \ge 0$.

\section*{Appendix C}

We first focus on  the analysis for any $i \in {\cal V}_1$ and $k \ge 2$. Denote $ {{{[A_1^k]}_{is}}}$ as the entry of the matrix $A_1^k$. By using (\ref{1-1}) recursively, we obtain $x_i^k = \sum\nolimits_{s \in {{\cal V}_1}} $ $ {{{[A_1^k]}_{is}}} x_s^0 + \sum\nolimits_{r = 1}^{k - 1} {\sum\nolimits_{s \in {{\cal V}_1}} {{{[A_1^{k - r}]}_{is}}} } \varepsilon _{1s}^{r - 1} + \varepsilon _{1s}^{k - 1}$. Moreover, from (\ref{4.3}) and (\ref{1-1}), one has ${{\bar x}^{k}} = \sum\nolimits_{i \in {{\cal V}_1}} {\rho _{1i}}(\sum\nolimits_{s \in {{\cal V}_1}} {a_{1is}}{x_s^{k-1} + \varepsilon _{1i}^{k-1}} ) = $ ${{\bar x}^{k-1}} + \sum\nolimits_{i \in {{\cal V}_1}} {\rho _{1i}} {\varepsilon _{1i}^{k-1}} $ for all $k \ge 1$. By using this fact recursively, one obtains ${{\bar x}^k} = \sum\nolimits_{s \in {{\cal V}_1}} {{\rho _{1s}}} x_s^0 + \sum\nolimits_{r = 1}^{k - 1} \sum\nolimits_{s \in {{\cal V}_1}} $ ${{\rho _{1s}}}  \varepsilon _{1s}^{r - 1} + \sum\nolimits_{s \in {{\cal V}_1}} {{\rho _{1s}}\varepsilon _{1s}^{k - 1}}$ for all $k \ge 2$. It then follows that $\| {x_i^k - {{\bar x}^k}} \| \le \sum\nolimits_{s \in {{\cal V}_1}} {|{{[A_1^k]}_{is}} - {\rho _{1s}}|\| {x_s^0} \|} + \sum\nolimits_{r = 1}^{k - 1}\sum\nolimits_{s \in {{\cal V}_1}} $ $ {{|{{[A_1^{k - r}]}_{is}} - {\rho _{1s}}|} } \| {\varepsilon _{1s}^{r - 1}} \| + (1 - {\rho _{1i}})\| {\varepsilon _{1i}^{k - 1}} \| + \sum\nolimits_{s \in {{\cal V}_1},s \ne i} $ ${{\rho _{1s}}\| {\varepsilon _{1s}^{k - 1}} \|}$, $\forall $$k \ge 2$. According to Proposition 1 in \cite{MaiVS16}, we know that there exist $\beta >0$ and $0<\gamma <1$ such that for any $i, s \in {\cal V}_1$ and $k \ge 0$, $\left| {{{[A_1^k]}_{is}} - {\rho _{1s}}} \right| \le \beta {\gamma ^k}$. Using this fact for the first two terms on the right-hand side of the above inequality and the fact $0 < {\rho _{1s}} < 1$, $\forall s \in {{\cal V}_1} $, for the last two terms, one has that for any $k \ge 2$,
\begin{equation}
\left\| {x_i^k - {{\bar x}^k}} \right\| \le \mathchar'26\mkern-10mu\lambda _1^k + \mathchar'26\mkern-10mu\lambda _2^k + \mathchar'26\mkern-10mu\lambda _3^k, \label{35}
\end{equation}
where $\mathchar'26\mkern-10mu\lambda _1^k = \beta {\gamma ^k}\sum\nolimits_{s \in {{\cal V}_1}} {\left\| {x_s^0} \right\|}$, $\mathchar'26\mkern-10mu\lambda _2^k = \beta \sum\nolimits_{r = 1}^{k - 1} {\gamma ^{k - r}}\sum\nolimits_{s \in {{\cal V}_1}} $ ${\left\| {\varepsilon _{1s}^{r - 1}} \right\|} $,
$\mathchar'26\mkern-10mu\lambda _3^k = \sum\nolimits_{s \in {{\cal V}_1}} \left\| {\varepsilon _{1s}^{k - 1}} \right\|$. By taking the square and the total expectation on both sides of (\ref{35}), one obtains $\Xi[{\left\| {x_i^k - {{\bar x}^k}} \right\|^2}] \le 3{(\mathchar'26\mkern-10mu\lambda _1^k)^2} + \Xi[3{(\mathchar'26\mkern-10mu\lambda _2^k)^2}] +  \Xi[3{(\mathchar'26\mkern-10mu\lambda _3^k)^2}]$, $\forall $$k \ge 2$, where the fact ${(\mathchar'26\mkern-10mu\lambda _1^k+\mathchar'26\mkern-10mu\lambda _2^k+\mathchar'26\mkern-10mu\lambda _3^k)^2}  \le 3[{(\mathchar'26\mkern-10mu\lambda _1^k)^2} + {(\mathchar'26\mkern-10mu\lambda _2^k)^2} + {(\mathchar'26\mkern-10mu\lambda _3^k)^2}]$ has been used. We next show that  $\sum\nolimits_{k = 2}^\infty  \Xi[{\left\| {x_i^k - {{\bar x}^k}} \right\|^2}] < \infty $, $\forall i \in {\cal V}_1$.

%?????????The proof of this theore can be reduced by replacing this by the previous latex that I have edited bofore
%
%??????This theorem is also in expectation, and all words refer to the almost sure should be changed into in expectation or learn how Hu Guoqiang describe this if it is not use in expectaion
%
%???????see gradient-free-6

Given an arbitrary $\ell \ge 2$. It is not hard to show that $\sum\nolimits_{k = 2}^{\ell  } 3{(\mathchar'26\mkern-10mu\lambda _1^k)^2} =3{\beta ^2}(\sum\nolimits_{s \in {{\cal V}_1}}{ {\left\| {x_s^0} \right\|} )^2}\frac{{{{\gamma ^4}}}(1 - {\gamma ^{2(\ell  - 1)}})}{{{1 - {\gamma ^2}}}}$. Note that $3{(\mathchar'26\mkern-10mu\lambda _2^k)^2} = 3{\beta ^2}\sum\nolimits_{r = 1}^{k - 1}\sum\nolimits_{r' = 1}^{k - 1}{{{\gamma ^{2k - r - r'}}\sum\nolimits_{s \in {{\cal V}_1}} {\sum\nolimits_{s' \in {{\cal V}_1}} {\| {\varepsilon _{1s}^{r - 1}} \|} } } } $ $ \cdot \| {\varepsilon _{1s'}^{r' - 1}} \|$. It follows from the relation $\| {\varepsilon _{1s}^{r - 1}} \|\big\| {\varepsilon _{1s'}^{r' - 1}} \big\| \le $ $ \frac{1}{2} ({\| {\varepsilon _{1s}^{r - 1}} \|^2}+ {\big\| {\varepsilon _{1s'}^{r' - 1}} \big \|^2})$ that $3{(\mathchar'26\mkern-10mu\lambda _2^k)^2} \le 3{m_1}{\beta ^2}\sum\nolimits_{r = 1}^{k - 1} \sum\nolimits_{r' = 1}^{k - 1}$ ${\gamma ^{2k - r - r'}} \sum\nolimits_{s \in {{\cal V}_1}} {{{\| {\varepsilon _{1s}^{r - 1}} \|}^2}} $. By virtue of Lemma {\ref{lem4}}  and Lemma {\ref{lem1}} (iv), one has $\Xi [{\left\| {\varepsilon _{1s}^{r - 1}} \right\|^2}] \le \frac{{\alpha _{r - 1}^2}}{{{{(\upsilon _{ss}^{r - 1})}^2}}}M_1^2$, $\forall r \ge 1$. It follows from both Lemma 1 and the proof of Proposition 1 in \cite{MaiVS16} that $\frac{1}{{\upsilon _{ss}^r}}>0$ is bounded for all $s \in {\cal V}_1$ and $r \ge 0$. In other words,  there exists a finite positive constant $\tilde \upsilon$ such that $ \frac{1}{{{{\upsilon _{ss}^{r-1}}}}} \le \tilde \upsilon$,  $\forall r \ge 1$. Thus, we have $\Xi [{\left\| {\varepsilon _{1s}^{r - 1}} \right\|^2}]  \le \alpha _{r - 1}^2{{\tilde \upsilon }^2}M_1^2$, $r \ge 1$. It then follows that the total expectation of $3{(\mathchar'26\mkern-10mu\lambda _2^k)^2}$ satisfies that $ \Xi [3{(\mathchar'26\mkern-10mu\lambda _2^k)^2}] \le \hbar \sum\nolimits_{r = 1}^{k - 1} {\alpha _{r - 1}^2\sum\nolimits_{r' = 1}^{k - 1} {{\gamma ^{2k - r - r'}}} } = \frac{\hbar }{{1 - \gamma }}\sum\nolimits_{r = 1}^{k - 1}$ $ {\alpha _{r - 1}^2} ( {\gamma ^{k - r + 1}}-{\gamma ^{2k - r}})$, where $ \hbar = 3m_1^2{\beta ^2}{{\tilde \upsilon }^2}M_1^2$. We thus know that $\sum\nolimits_{k = 2}^{\ell } {\Xi [3{{(\lambda _2^k)}^2}]}  \le \frac{\hbar }{{1 - \gamma }}\sum\nolimits_{k = 2}^{\ell } \sum\nolimits_{r = 1}^{k - 1} {\alpha _{r - 1}^2} ({\gamma ^{k - r + 1}}-{\gamma ^{2k - r}} ) $. Note that $\sum\nolimits_{k = 2}^{\ell  } \sum\nolimits_{r = 1}^{k - 1}{\alpha _{r - 1}^2}({\gamma ^{k - r + 1}}- {\gamma ^{2k - r}}) = \sum\nolimits_{r = 1}^{\ell  - 1} \alpha _{r - 1}^2$ $\sum\nolimits_{k = r + 1}^{\ell }({\gamma ^{k - r + 1}} - {\gamma ^{2k - r}})   = \sum\nolimits_{r = 1}^{\ell  - 1}\frac{\alpha _{r - 1}^2}{1-{\gamma ^2}}({\gamma ^{2\ell  - r+2}}-{\gamma ^{\ell  - r + 3}} $ ${{ - {\gamma ^{\ell - r + 2 }} - {\gamma ^{r + 2}} + {\gamma ^2} + {\gamma ^3}}}) \le \frac{{{1 + {\gamma ^2} + {\gamma ^3}}}}{{{1 - {\gamma ^2}}}} \sum\nolimits_{r = 1}^{\ell  - 1} {\alpha _{r - 1}^2}$, where the inequality follows from the facts $0<\gamma <1$ and $ 2 \ell -r +2 >0$. Therefore, we have $\sum\nolimits_{k = 2}^{\ell } {\Xi [3{{(\lambda _2^k)}^2}]}  \le \frac{ \hbar (1 + {\gamma ^2} + {\gamma ^3})}{{{(1 - \gamma )(1 - {\gamma ^2})}}} \sum\nolimits_{r = 1}^{\ell  - 1} {\alpha _{r - 1}^2} $. Moreover, by noting that $ 3{(\mathchar'26\mkern-10mu\lambda _3^k)^2} \le $ $3{m_1}{\sum\nolimits_{s \in {{\cal V}_1}} {\| {\varepsilon _{1s}^{k - 1}} \|} ^2}$, it can be shown that $\sum\nolimits_{k = 2}^{\ell }{\Xi [3{{(\mathchar'26\mkern-10mu\lambda _3^k)}^2}]}  \le 3m_1^2{\tilde \upsilon ^2}M_1^2\sum\nolimits_{k = 2}^{\ell } {\alpha _{k - 1}^2}$. Based on the above analysis, we obtain $\sum\nolimits_{k = 2}^{\infty}\{3{(\mathchar'26\mkern-10mu\lambda _1^k)^2}+\Xi[3{(\mathchar'26\mkern-10mu\lambda _2^k)^2}] +  \Xi[3{(\mathchar'26\mkern-10mu\lambda _3^k)^2}] \} < \infty $ as $\ell  \to \infty $, where we have used the facts $0< \gamma <1 $ and $\sum\nolimits_{k = 0}^\infty  {\alpha _k^2} < \infty $. Therefore, we have $\sum\nolimits_{k = 2}^\infty  \Xi[{\left\| {x_i^k - {{\bar x}^k}} \right\|^2}] < \infty$,  $\forall i \in {\cal V}_1$.

Following the similar analysis to obtain  (\ref{35}), we can derive that for any $ i \in {\cal V}_1$, $\left\| {x_i^1 - {{\bar x}^1}} \right\| \le \mathchar'26\mkern-10mu\lambda _1^1 +\mathchar'26\mkern-10mu\lambda _3^1 ,$
where $\mathchar'26\mkern-10mu\lambda _1^1 = \beta {\gamma }\sum\nolimits_{s \in {{\cal V}_1}} {\left\| {x_s^0} \right\|}$,
$\mathchar'26\mkern-10mu\lambda _3^1 = \sum\nolimits_{s \in {{\cal V}_1}} \left\| {\varepsilon _{1s}^{0}} \right\|$. Based on this fact, it can be shown that $\Xi [{\left\| {x_i^1 - {{\bar x}^1}} \right\|^2}] < \infty$,  $\forall i \in {\cal V}_1$. Summarizing all the above analysis, we have $\sum\nolimits_{k = 0}^\infty  {\Xi [{{\left\| {x_i^k - {{\bar x}^k}} \right\|}^2}]}  < \infty $,  $\forall i \in {\cal V}_1$. This also implies that given any $i \in {\cal V}_1$, one has ${\Xi [{{\left\| {x_i^k - {{\bar x}^k}} \right\|}^2}]}  < \infty $ for all $k \ge 0$. Based on these two facts and Theorem 4.2.1 in \cite{Lukacs14}, we know that for each  $i \in {\cal V}_1$, the sequence $\{\sum\nolimits_{k = 0}^\ell  {{\left\| {x_i^k - {{\bar x}^k}} \right\|}^2} \}_{\ell \ge 0}$ is almost surely convergent. This implies that  for each  $ i \in {\cal V}_1$, the limit of ${\left\| {x_i^k - {{\bar x}^k}} \right\|}^2$ exists  with probability one as $k \to \infty$, and ${\lim _{k \to \infty }}{\left\| {x_i^k - {{\bar x}^k}} \right\|}=0$ with probability one, i.e., all the agents in ${\cal V}_1$ achieve consensus almost surely.

% Furthermore, we have $\Xi [{\left\| {x_i^1 - {{\bar x}^1}} \right\|^2}] \le 2{(\mathchar'26\mkern-10mu\lambda _1^1)^2} + \Xi [2{(\mathchar'26\mkern-10mu\lambda _3^1)^2}]$. Using the similar analysis to obtain the upper bound of ?????????$\Xi [3{(\mathchar'26\mkern-10mu\lambda _3^k)^2}]$ for $k \ge 2$, we can derive the upper bound of  $\Xi [2{(\mathchar'26\mkern-10mu\lambda _3^1)^2}]$ which is  shown to be finite.
%Then, it is not hard to show that $\Xi [{\left\| {x_i^1 - {{\bar x}^1}} \right\|^2}] < \infty$,  $\forall i \in {\cal V}_1$.

Applying the similar analysis above for $ i \in {\cal V}_1$ to $j  \in {\cal V}_2$, one can obtain $\sum\nolimits_{k = 0}^\infty  {\Xi [{{\left\| {y_j^k - {{\bar y}^k}} \right\|}^2}]}  < \infty $,  $\forall j \in {\cal V}_2$. Furthermore, we can derive that for each ${j \in {{\cal V}_2}}$, $\mathop {\lim }\nolimits_{k \to \infty } \| y_j^k - {\bar y^k} \|  = 0$ with probability one,  i.e., all the agents in ${\cal V}_2$ achieve consensus almost surely.

\section*{Appendix D}

For the simplicity of the notation, we will use $g_{1i\mu ,k}$ to represent $g_{1i\mu}(x_i^k,\Pi _{1i}^k)$ for the remaining proofs. In the proof of Theorem {\ref{thm1}}, one obtained that for any $ k \ge 0$, ${{\bar x}^{k + 1}} = {{\bar x}^k}  + \sum\nolimits_{i \in {{\cal V}_1}} {{\rho _{1i}}\varepsilon _{1i}^k} $. Then, it follows that for any $(x^*_{\mu},y^*_{\mu}) \in {\cal Z}^*_{\mu}$ and $k \ge 0$, ${\left\| {{{\bar x}^{k + 1}} - x_{\mu}^*} \right\|^2}={\left\| {{{\bar x}^k} - x_{{\mu }}^*} \right\|^2} + \|\sum\nolimits_{i \in {{\cal V}_1}} {{\rho _{1i}}\varepsilon _{1i}^k} \|^2+ 2{({{\bar x}^k} - x_{\mu}^*)^T}(\sum\nolimits_{i \in {{\cal V}_1}} {{\rho _{1i}}\varepsilon _{1i}^k} )$. Taking the conditional expectation for both sides of this equality on ${\cal F}_k$ yields
\begin{eqnarray}
& &\Xi [{\left\| {{{\bar x}^{k + 1}} - x_{{\mu}}^*} \right\|^2}\big|{{\cal F}_k}] \notag \\
& =& {\left\| {{{\bar x}^k} - x_{{\mu }}^*} \right\|^2} + \Xi [\big\|\sum\nolimits_{i \in {{\cal V}_1}} {{\rho _{1i}}\varepsilon _{1i}^k} \big\|^2\big|{{\cal F}_k}]  \notag \\
& &+ \Xi [ 2{({{\bar x}^k} - x_{{\mu}}^*)^T}(\sum\nolimits_{i \in {{\cal V}_1}} {{\rho _{1i}}\varepsilon _{1i}^k} )\big|{{\cal F}_k}]. \label{4.6}
\end{eqnarray}
Let ${\phi^k _1}= {\sum\nolimits_{i \in {{\cal V}_1}} 2{{\rho _{1i}}{({{\bar x}^k} - x_{\mu}^*)^T}(\varepsilon _{1i}^k + \frac{{{\alpha _k}}}{{\upsilon _{ii}^k}}g_{1i\mu,k})} }$, ${\phi^k _2}=$ $\sum\nolimits_{i \in {{\cal V}_1}} 2{{\alpha _k}{{({{\bar x}^k} - x_{\mu}^*)}^T}\frac{1}{{\upsilon _{ii}^k}}(\upsilon _{ii}^k - } {\rho _{1i}})g_{1i\mu,k}$, ${\phi^k _3} = - \sum\nolimits_{i \in {{\cal V}_1}} $ $2{\alpha _k}{{({{\bar x}^k} - x_{\mu}^*)^T}g_{1i\mu,k}}$. It is clear that $2{({{\bar x}^k} - x_{{\mu}}^*)^T}(\sum\nolimits_{i \in {{\cal V}_1}}  $ ${{\rho _{1i}}\varepsilon _{1i}^k} )$ $={\phi^k _1}+{\phi^k _2}+{\phi^k _3}$.

For ${\phi^k _1}$, we have ${\phi^k _1} = 2\sum\nolimits_{i \in {{\cal V}_1}} {{\rho _{1i}}} ({{\bar x}^k} - {{\bar x}^{k + 1}}+{{\bar x}^{k + 1}} $ ${- x_i^{k + 1}+x_i^{k + 1} - x_{\mu}^*)^T}(\varepsilon _{1i}^k + \frac{{{\alpha _k}}}{{\upsilon _{ii}^k}}g_{1i\mu,k}) \le 2\sum\nolimits_{i \in {{\cal V}_1}}{\rho _{1i}}$ $ \cdot {(\left\| {{{\bar x}^k} - {{\bar x}^{k + 1}}} \right\|+\left\| {{{\bar x}^{k + 1}} - x_i^{k + 1}} \right\|}) \big\|\varepsilon _{1i}^k + \frac{{{\alpha _k}}}{{\upsilon _{ii}^k}}g_{1i\mu,k}\big\|$,
where  the fact that ${(x_i^{k + 1} - x_{\mu}^*)^T}(\varepsilon _{1i}^k + \frac{{{\alpha _k}}}{{\upsilon _{ii}^k}}g_{1i\mu,k}) \le 0$ based on Lemma 1 (a) in \cite{Nedic10} has been used in the inequality. Then,
\begin{eqnarray}
\Xi [\left. {\phi _1^k} \right|{{\cal F}_k}] \le  \zeta^k _{1} + \zeta^k _2,   \quad \label{10}
\end{eqnarray}
where $\zeta^k _1=\Xi [4\sum\nolimits_{i \in {{\cal V}_1}} {{\rho _{1i}}} \| {{{\bar x}^k} - {{\bar x}^{k + 1}}} \|\frac{{{\alpha _k}}}{{\upsilon _{ii}^k}}\left\| {g_{1i\mu,k}} \right\||{{\cal F}_k}] $, $\zeta^k _2$ $=\Xi [4\sum\nolimits_{i \in {{\cal V}_1}} {{\rho _{1i}}} \big\| {{{\bar x}^{k + 1}} - x_i^{k + 1}} \big\|\frac{{{\alpha _k}}}{{\upsilon _{ii}^k}}\big\| {g_{1i\mu,k}} \big\||{{\cal F}_k}]$, with Lemma {\ref{lem4}} being used.

For $\phi _2^k$, it follows that  $\phi _2^k \le \sum\nolimits_{i \in {{\cal V}_1}} {[\frac{{\alpha _k^2}}{{{m_1}}}{{\left\| {{{\bar x}^k} - x_\mu ^*} \right\|}^2}}+ {m_1}{{\tilde \upsilon }^2}$ $ \cdot | \upsilon _{ii}^k - {\rho _{1i}} |^2{\left\| {{g_{1i\mu ,k}}} \right\|^2}] \le \alpha _k^2{\left\| {{{\bar x}^k} - x_\mu ^*} \right\|^2} + {m_1}{{\tilde \upsilon }^2}$ ${\beta ^2}{\gamma ^{2k}}\sum\nolimits_{i \in {{\cal V}_1}}$ ${\left\| {{g_{1i\mu ,k}}} \right\|} ^2$, where the relation $| {\upsilon _{ii}^k - {\rho _{1i}}} | \le \beta {\gamma ^k}$ (see Proposition 1 in \cite{MaiVS16}) has been used in the second inequality, and $\beta$ and $\gamma $ have been given in the proof of Theorem {\ref{thm1}}. Taking the conditional expectation for both sides of this equation on ${\cal F}_k$ and using Lemma {\ref{lem1}} (iv) yield
\begin{equation}
\Xi [\left.\phi _2^k\right|{{\cal F}_k}] \le \alpha _k^2{\left\| {{{\bar x}^k} - x_{\mu}^*} \right\|^2} + m_1^2{{\tilde \upsilon } ^2}{\beta ^2}M_1^2{\gamma ^{2k}}.  \label{11}
\end{equation}

For ${\phi^k _3}$, it follows that  $\phi _3^k =  - \sum\nolimits_{i \in {{\cal V}_1}} 2{\alpha _k}({{\bar x}^k} - x_i^k+x_i^k -$ $ x_{\mu}^*)^Tg_{1i\mu,k}$. Let $\zeta^k _{3}=M_1\sum\nolimits_{i \in {{\cal V}_1}}2{\alpha _k} \left\| {{{\bar x}^k} - x_i^k} \right\|$. By taking the conditional expectation for $\phi _3^k$ on ${\cal F}_k$, one obtains $ \Xi [\phi _3^k|{{\cal F}_k}] =  - \sum\nolimits_{i \in {{\cal V}_1}} 2 {\alpha _k}{({{\bar x}^k} - x_i^k)^T}\Xi [{g_{1i\mu ,k}}|{{\cal F}_k}] - \sum\nolimits_{i \in {{\cal V}_1}} $ $2{\alpha _k} {{{(x_i^k - x_\mu ^*)}^T}}\Xi [{g_{1i\mu ,k}}|{{\cal F}_k}]\le \zeta^k _{3} +2{\alpha _k} \sum\nolimits_{i \in {{\cal V}_1}}[ f_{1i{\mu}}(x_{\mu}^*, \Pi _{1i}^k)$ $- f_{1i{\mu}}(x_i^k,\Pi _{1i}^k) ] $, where  the relations $\Xi [g_{1i\mu,k}\big|{{\cal F}_k} ] ={\nabla _{x_i^k}}f_{1i\mu}(x_i^k,\Pi _{1i}^k)$ and $\| {{\nabla _{x_i^k}}{f_{1i\mu }}(x_i^k,\Pi _{1i}^k)} \| \le M_1$ and the convexity of ${f_{1i\mu}}( \cdot ,\Pi _{1i}^k)$ at $x_i^k$, which are based on Lemma {\ref{lem1}} (i)-(iii), have been used in the inequality.

To present the upper bound of $2{\alpha _k} \sum\nolimits_{i \in {{\cal V}_1}}[ f_{1i{\mu}}(x_{\mu}^*,\Pi _{1i}^k)- f_{1i{\mu}}(x_i^k,\Pi _{1i}^k) ]$, we introduce ${\tilde \zeta^k _{4}} = 2{\alpha _k} \sum\nolimits_{i \in {{\cal V}_1}}[f_{1i{\mu}}(x_{\mu}^*,\Pi _{1i}^k)- f_{1i{\mu}}(x_{\mu}^*, {{\bar y}^k})]$, ${\tilde \zeta^k _{5}} =2{\alpha _k}\sum\nolimits_{i \in {{\cal V}_1}} [f_{1i{\mu}}({{\bar x}^k},{{\bar y}^k}) -f_{1i{\mu}}(x_i^k,{{\bar y}^k})]$, and ${\tilde \zeta^k _{6}} = 2{\alpha _k}\sum\nolimits_{i \in {{\cal V}_1}} [f_{1i{\mu}}(x_i^k,{{\bar y}^k}) - f_{1i{\mu}}(x_i^k,\Pi _{1i}^k)]$. It is clear that $2{\alpha _k} \sum\nolimits_{i \in {{\cal V}_1}}[ f_{1i{\mu}}(x_{\mu}^*,\Pi _{1i}^k)- f_{1i{\mu}}(x_i^k,\Pi _{1i}^k) ]= {\tilde \zeta^k _{4}}+ 2{\alpha _k}\sum\nolimits_{i \in {{\cal V}_1}} [f_{1i{\mu}}(x_{\mu}^*,{{\bar y}^k}) -f_{1i\mu}({{\bar x}^k},{{\bar y}^k})]+{\tilde \zeta^k _{5}}+{\tilde \zeta^k _{6}}$. For ${\tilde \zeta^k _{4}}$, it follows that ${\tilde \zeta^k _{4}} \le 2{\alpha _k}\sum\nolimits_{i \in {{\cal V}_1}}{\nabla _{{{\bar y}^k}}}f_{1i{\mu}}(x_{\mu}^*,{{\bar y}^k})(\Pi _{1i}^k - {{\bar y}^k}) \le 2{\alpha _k}\sum\nolimits_{i \in {{\cal V}_1}} \left\| {{\nabla _{{{\bar y}^k}}}f_{1i{\mu}}(x_{\mu}^*,{{\bar y}^k})} \right\|\big\| \frac{1}{{\sum\nolimits_{s \in {{\cal V}_2}} {{a_{1is}}} }}\sum\nolimits_{s \in {{\cal V}_2}} {{a_{1is}}y_s^k} - {{\bar y}^k}\big\|=2{\alpha _k}\sum\nolimits_{i \in {{\cal V}_1}} \left\| {{\nabla _{{{\bar y}^k}}}f_{1i{\mu}}(x_{\mu}^*,{{\bar y}^k})} \right\| \big\|\sum\nolimits_{s \in {{\cal V}_2}} \frac{{{a_{1is}}}}{{\sum\nolimits_{s \in {{\cal V}_2}} {{a_{1is}}} }}$ $(y_s^k  - {{\bar y}^k}) \big\|\le \zeta ^k _{4}$, where we have used the concavity of $f_{1i{\mu}}(x_{\mu}^*,{{\cdot}})$ at ${\bar y}^k$ in the first inequality, the definition of $\Pi _{1i}^k$ in the second inequality, and Lemma {\ref{lem1}} (iii) in the last inequality, with $\zeta ^k _{4}=2{\alpha _k}m_1{M_2}{\sum\nolimits_{s \in {{\cal V}_2}} \left\| {y_s^k - {{\bar y}^k}} \right\|}$. Similarly, it can be shown that ${\tilde \zeta^k _{5}} \le 2{\alpha _k}\sum\nolimits_{i \in {{\cal V}_1}} {\nabla _{{{\bar x}^k}}}f_{1i\mu}({{\bar x}^k},{{\bar y}^k})({{\bar x}^k}-$ $ x_i^k) \le \zeta ^k _{5}$ and ${\tilde \zeta^k _{6}} \le 2{\alpha _k}\sum\nolimits_{i \in {{\cal V}_1}} {{\nabla _{\Pi _{1i}^k}}f_{1i{\mu}}(x_i^k,\Pi _{1i}^k)({{\bar y}^k} - \Pi _{1i}^k)} $ $\le \zeta ^k _{4}$, where $\zeta ^k _{5}=2{\alpha _k}{M_1}\sum\nolimits_{i \in {{\cal V}_1}}\| {{\bar x}^k} -{ x_i^k\|}$. It then follows that $2{\alpha _k} \sum\nolimits_{i \in {{\cal V}_1}}[ f_{1i{\mu}}(x_{\mu}^*,\Pi _{1i}^k)- f_{1i{\mu}}(x_i^k,\Pi _{1i}^k) ] \le 2{\zeta^k _{4}} +{\zeta^k _{5}} +$ $2{\alpha _k}\sum\nolimits_{i \in {{\cal V}_1}} [f_{1i{\mu}}(x_{\mu}^*,{{\bar y}^k}) -f_{1i\mu}({{\bar x}^k},{{\bar y}^k})]$. Therefore,
\begin{align}
& \Xi[{\phi^k _3} \big| {\cal F}_k] \le  \zeta^k _{3} +2{\zeta^k _{4}} +{\zeta^k _{5}} \notag \\
& + 2{\alpha _k}\sum\nolimits_{i \in {{\cal V}_1}} [f_{1i{\mu}}(x_{\mu}^*,{{\bar y}^k}) - f_{1i{\mu}}({{\bar x}^k},{{\bar y}^k})]. \label{22}
\end{align}

Note that $\Xi [{\left\| \sum\nolimits_{i \in {{\cal V}_1}} {{\rho _{1i}}}{\varepsilon _{1i}^k} \right\|^2}\big|{{\cal F}_k}] \le \Xi [\sum\nolimits_{i \in {{\cal V}_1}} {{\rho _{1i}}}\|{\varepsilon _{1i}^k \|^2} $ $\big|{{\cal F}_k}]  \le \Xi [\sum\nolimits_{i \in {{\cal V}_1}} {{\rho _{1i}}}\frac{{\alpha _k^2}}{{{{(\upsilon _{ii}^k)}^2}}}{\| {{g_{1i\mu ,k}}}\|^2}|{{\cal F}_k}] \le \alpha _k^2{\tilde \upsilon ^2}M_1^2$, where we have used the fact that $\left\|  \cdot  \right\|^2$ is a convex function  in the first inequality, Lemma {\ref{lem4}} in the second inequality, and  Lemma {\ref{lem1}} (iv) and $\sum\nolimits_{i \in {{\cal V}_1}} {{\rho _{1i}}}=1$ in the third inequality.  Substituting the above inequality  and the equations (\ref{10})-(\ref{22}) into the right-hand side of (\ref{4.6}) yields that  for any $k \ge 0$,
\begin{align}
& \Xi [ {\left\| {{{\bar x}^{k + 1}} - x_{\mu}^*} \right\|^2} \big|{{\cal F}_k}] \le (1 + \alpha _k^2){\left\| {{{\bar x}^k} - x_{\mu}^*} \right\|^2}+p_1^k \notag \\
& +2{\alpha _k}[F_{{\mu}}(x_{\mu}^*,{{\bar y}^k}) - F_{\mu}({{\bar x}^k},{{\bar y}^k})], \label{15}
\end{align}
where  $p_1^k=2\sum\nolimits_{\iota  = 1}^5 {\varsigma _\iota ^k}+{\tilde \upsilon ^2}M_1^2\alpha _k^2+ m_1^2{{\tilde \upsilon } ^2}{\beta ^2}M_1^2{\gamma ^{2k}}$. Next, we show that $\sum\nolimits_{k = 0}^{\infty}\Xi [p^k_1]  < \infty $. We first consider $\sum\nolimits_{k = 0}^{\infty}\Xi [{2\varsigma _2^k}]$. For any $ k \ge 0$, one has $\zeta _2^k \le \Xi [2\sum\nolimits_{i \in {{\cal V}_1}} {{\rho _{1i}}} ({\left\| {{{\bar x}^{k + 1}} - x_i^{k + 1}} \right\|^2} + \frac{{\alpha _k^2}}{{{{(\upsilon _{ii}^k)}^2}}}{\left\| {{g_{1i\mu ,k}}} \right\|^2})|{{\cal F}_k}] \le 2\sum\nolimits_{i \in {{\cal V}_1}} {{\rho _{1i}}} \Xi [{\left\| {{{\bar x}^{k + 1}} - x_i^{k + 1}} \right\|^2}|{{\cal F}_k}] + 2\alpha _k^2{{\tilde \upsilon }^2}M_1^2$. Then, we know that $\sum\nolimits_{k = 0}^\infty  {\Xi [2\zeta _2^k]}  \le 4\sum\nolimits_{i \in {{\cal V}_1}} {{\rho _{1i}}} \sum\nolimits_{k = 0}^\infty  \Xi [\| {{\bar x}^{k + 1}} - $ $x_i^{k + 1} \|^2]  + 4{{\tilde \upsilon }^2}M_1^2\sum\nolimits_{k = 0}^\infty  {\alpha _k^2}$. It has been shown in the proof of Theorem {\ref{thm1}} that $\sum\nolimits_{k = 0}^\infty \Xi [{\left\| {{{\bar x}^{k}} - x_i^{k}} \right\|^2}] < \infty$ for each $ i \in {\cal V}_1$. Thus, we know that $\sum\nolimits_{k = 0}^\infty  {\Xi [2 \zeta _2^k]}  < \infty $. For $\iota  = 1,3,4,5$, we can similarly derive $\sum\nolimits_{k = 0}^\infty  \Xi  [2\varsigma _\iota ^k] \le \hat \varsigma _\iota ^k + \hat {\hbar} _\iota \sum\nolimits_{k = 0}^\infty  {\alpha _k^2}$, where $\hat \varsigma _\iota ^k $ satisfies $\hat \varsigma _\iota ^k  < \infty $, $\hat {\hbar} _\iota $ is a finite positive constant. This implies that $\sum\nolimits_{k = 0}^\infty  \Xi  [2\varsigma _\iota ^k] < \infty $, where $\iota = 1,3,4,5$. Moreover, we know that $\sum\nolimits_{k = 0}^\infty {\gamma ^{2k}} \le \frac{1}{{1 - {\gamma ^2}}} < \infty $, and therefore $\sum\nolimits_{k = 0}^{\infty}\Xi [p^k_1]  < \infty $.

On the other hand, we know that ${{\bar y}^{k + 1}} = {{\bar y}^k}  + \sum\nolimits_{j \in {{\cal V}_2}} {{\rho _{2j}}\varepsilon _{2j}^k} $, $\forall k \ge 0$. It follows from the similar analysis to obtain (\ref{15}) and $\sum\nolimits_{k = 0}^{\infty} \Xi [p^k_1] < \infty $  that  for any $k \ge 0$,
\begin{align}
& \Xi [ {\left\| {{{\bar y}^{k + 1}} - y_{\mu}^*} \right\|^2}\big|{{\cal F}_k}] \le (1 + \alpha _k^2){\left\| {{{\bar y}^k} - y_{\mu}^*} \right\|^2} + p_2^k \nonumber \\
& - 2{\alpha _k}[{{F_{\mu}}({{\bar x}^k},y_{\mu}^*) - {F_{\mu}}({{\bar x}^k},{{\bar y}^k})}], \label{23}
\end{align}
where $ p^k_2$ is nonnegative and satisfying $\sum\nolimits_{k = 0}^{\infty}$ $\Xi [p^k_2] < \infty $.

Summing up (\ref{15}) and (\ref{23}) yields that  for any $k \ge 0$,
\begin{eqnarray}
& &\Xi [ ({\left\| {{{\bar x}^{k + 1}} - x_{\mu}^*} \right\|^2}+{\left\| {{{\bar y}^{k + 1}} - y_{\mu}^*} \right\|^2})  \big|{{\cal F}_k}]  \notag \\
&\le& (1 + \alpha _k^2)({\left\| {{{\bar x}^k} - x_{\mu}^*} \right\|^2+\left\| {{{\bar y}^k} - y_{\mu}^*} \right\|^2}) + p^k-q^k,  \qquad \ \label{36}
\end{eqnarray}
where $p^k=p_1^k+p_2^k$, $q^k=2{\alpha _k}[ -F_{\mu}(x_{\mu}^*,{{\bar y}^k}) + F_{\mu}(x_{\mu}^*,y_{\mu}^*) + F_{\mu}({{\bar x}^k},$ $y_{\mu}^*)- F_{\mu}(x_{\mu}^*,y_{\mu}^*)]$. It is obvious that $\sum\nolimits_{k = 0}^\infty  {\Xi [{p^k}]} < \infty$.
\medskip

\ifCLASSOPTIONcaptionsoff
  \newpage
\fi

\end{document}